\documentclass[a4paper,12pt,reqno]{amsart}
\usepackage{amsmath}
\usepackage{amssymb}
\usepackage{amsthm}

\usepackage{amscd}
\usepackage{amsfonts}
\usepackage{setspace}
\usepackage{version}

\usepackage{hyperref}

\title[Lower bounds for low moments of character sums, I]{Lower bounds for low moments of character sums, I: Short sums with general multiplicative weights}
\author{Adam J Harper}
\address{Mathematics Institute, Zeeman Building, University of Warwick, Coventry CV4 7AL, England}
\email{A.Harper@warwick.ac.uk}
\date{1st July 2026}
\thanks{This research was funded in part by the Engineering and Physical Sciences Research Council of the United Kingdom [grant EP/V055755/1]. Support was also received from the Swedish Research Council [grant no. 2021-06594], while the author was in residence at Institut Mittag-Leffler, Djursholm, Sweden, during the 2024 Analytic Number Theory program; and from the Simons Foundation and the Centre de Recherches Math\'ematiques, Montr\'eal, while the author was in residence as Aisenstadt Chair during the 2026 Universal Statistics in Number Theory thematic semester. For the purpose of open access, the author has applied a Creative Commons Attribution (CC-BY) licence to any Author Accepted Manuscript version arising from this submission.}

\numberwithin{equation}{section}
\theoremstyle{plain}

\onehalfspace
\newcommand{\N}{\mathbb{N}}
\newcommand{\R}{\mathbb{R}}
\newcommand{\E}{\mathbb{E}}

\newcommand{\Echar}{\mathbb{E}^{\text{char}}}

\newtheorem{thmlb1}{Theorem}
\newtheorem{thmlb2}[thmlb1]{Theorem}
\newtheorem{thmlb3}[thmlb1]{Theorem}

\newtheorem{keyprop1}{Key Proposition}
\newtheorem{keyprop2}[keyprop1]{Key Proposition}

\newtheorem{approxres1}{Approximation Result}

\newtheorem{numth1}{Number Theory Result}
\newtheorem{numth2}[numth1]{Number Theory Result}
\newtheorem{numth3}[numth1]{Number Theory Result}
\newtheorem{harman1}{Harmonic Analysis Result}

\newtheorem{probres1}{Probability Result}
\newtheorem{probres2}[probres1]{Probability Result}
\newtheorem{probres3}[probres1]{Probability Result}
\newtheorem{probres4}[probres1]{Probability Result}
\newtheorem{probres5}[probres1]{Probability Result}

\newtheorem{keyprop3}[keyprop1]{Key Proposition}
\newtheorem{keyprop4}[keyprop1]{Key Proposition}

\newtheorem{lem1}{Lemma}

\begin{document}

\maketitle

\begin{abstract}
We establish sharp lower bounds for the Dirichlet character moments $\frac{1}{r-1} \sum_{\chi \; \text{mod} \; r} |\sum_{n \leq x} \chi(n)|^{2q}$, where $r$ is a large prime, $1 \leq x \leq r^{0.499}$, and $0 \leq q \leq 1$ is real. These match the better than squareroot cancellation upper bounds obtained in previous work of the author. We prove the same sharp lower bounds for the moments $\frac{1}{T} \int_{0}^{T} |\sum_{n \leq x} n^{it}|^{2q} dt$ of zeta sums, and more generally for moments of character sums $\sum_{n \leq x} h(n) \chi(n)$ with suitably bounded multiplicative twist $h(n)$.

The proofs are based on a comparison of the sizes of $\frac{1}{r-1} \sum_{\chi \; \text{mod} \; r} (\sum_{n \leq x} \chi(n)) \overline{I(\chi)}$, $\frac{1}{r-1} \sum_{\chi \; \text{mod} \; r} |I(\chi)|^2$ and $\frac{1}{r-1} \sum_{\chi \; \text{mod} \; r} |I(\chi)|^4$, where $I(\chi)$ is a certain ``barrier adjusted'' Perron integral inspired by the analogous results for random multiplicative functions.

In a companion paper, we extend these arguments to the full interesting range $x \leq 0.99r$ for the unweighted character sum moments $\frac{1}{r-1} \sum_{\chi \; \text{mod} \; r} |\sum_{n \leq x} \chi(n)|^{2q}$. This leads to a positive proportion non-vanishing result for Dirichlet theta functions $\theta(1,\chi)$. 
\end{abstract}

\section{Introduction}
This paper is concerned with the average size, in the sense of moments, of various (possibly weighted) sums of Dirichlet characters $\chi(n)$ or ``continuous characters'' $n^{it}$. This is an extremely classical and well-studied problem, for example Montgomery and Vaughan~\cite{mvchar} proved that for any large $r \in \N$, any real $q > 0$ and any $x$, we have
$$ \frac{1}{\varphi(r)} \sum_{\chi \neq \chi_0 \; \text{mod} \; r} |\sum_{n \leq x} \chi(n)|^{2q} \ll_{q} r^q . $$
(In fact, they proved this with $\sum_{n \leq x} \chi(n)$ replaced by $M(\chi) := \max_{x} |\sum_{n \leq x} \chi(n)|$.) In some cases this bound will be sharp for all $q$, for example it would be if $r$ were prime and $0.01r \leq x \leq 0.99r$, say. But if $x$ is significantly smaller than $r$, the bound becomes less good. Granville and Soundararajan~\cite{gransoundlcs} investigated high integer moments, for various sizes of $x$, as a means of showing the existence of unusually large character sums. There are many other relevant papers, some others will be mentioned below but we do not attempt a full survey. We note, however, that Szab\'o~\cite{szaboupper,szabolower} has shown that for any large prime $r$, any real $q > 2$ and $1 \leq x \leq r/2$, we have
$$ \frac{1}{r-1} \sum_{\chi \neq \chi_0 \; \text{mod} \; r} |\sum_{n \leq x} \chi(n)|^{2q} \asymp_{q} x^q \log^{(q-1)^2}(10L) , $$
where $L = L_r := \min\{x, r/x\}$. (The upper bound part of this assumes the Generalised Riemann Hypothesis for Dirichlet $L$-functions, and is also established for non-prime $r$ if one replaces $\frac{1}{r-1}$ by $\frac{1}{\varphi(r)}$, and restricts the outer sum to the primitive characters mod $r$.) Thus we understand the order of magnitude, for moments of unweighted character sums $\sum_{n \leq x} \chi(n)$ higher than the fourth moment. For moments of character sums $\sum_{n \leq x} \mu(n) \chi(n)$ twisted by the M\"obius function $\mu(n)$, we also note Gorodetsky's paper~\cite{gorodetsky}, which makes precise conjectures about the integer moments based on function field and random matrix considerations.

Our focus shall be on {\em low power moments}, i.e. moments below the second moment. Letting $r$ be a large prime, the orthogonality of Dirichlet characters implies that for any $1 \leq x < r$, we have
$$ \frac{1}{r-1} \sum_{\chi \; \text{mod} \; r} |\sum_{n \leq x} \chi(n)|^{2} = \lfloor x \rfloor . $$
Using H\"older's inequality, it follows that $\frac{1}{r-1} \sum_{\chi \; \text{mod} \; r} |\sum_{n \leq x} \chi(n)|^{2q} \leq x^q$ for all $0 \leq q \leq 1$. This easy bound corresponds to squareroot cancellation on average, and on first inspection one might suppose it to be essentially sharp. But in previous work~\cite{harpertypicalchar}, the author showed that uniformly for $1 \leq x \leq r$ and $0 \leq q \leq 1$, we have
\begin{equation}\label{eqbtscupper}
\frac{1}{r-1} \sum_{\chi \; \text{mod} \; r} |\sum_{n \leq x} \chi(n)|^{2q} \ll \Biggl(\frac{x}{1 + (1-q)\sqrt{\log\log(10L)}} \Biggr)^q ,
\end{equation}
where again $L = \min\{x, r/x\}$. In particular, whenever $x \rightarrow \infty$ with $r$ but $x = o(r)$, we have $\frac{1}{r-1} \sum_{\chi \; \text{mod} \; r} |\sum_{n \leq x} \chi(n)|^{2q} = o(x^q )$ for any fixed $0 < q < 1$. Thus we have the unexpected phenomenon of {\em better than squareroot cancellation}.

This behaviour persists if one introduces a multiplicative twist into the sum, for example it was shown in \cite{harpertypicalchar} that the bound \eqref{eqbtscupper} is valid uniformly for all sums $\frac{1}{r-1} \sum_{\chi \; \text{mod} \; r} |\sum_{n \leq x} h(n) \chi(n)|^{2q}$, where $h(n)$ is any multiplicative function that has absolute value 1 on primes and absolute value at most 1 on prime powers. See also the recent paper of Gao and Wu~\cite{gaowulow}, which adapts the arguments of \cite{harpertypicalchar} to handle the case where $h(n)$ are the Fourier coefficients of a fixed Hecke eigenform. Likewise for continuous characters $n^{it}$, uniformly for any large $T$, any $1 \leq x \leq T$ and any $0 \leq q \leq 1$, we have
\begin{equation}\label{eqbtsccontupper}
\frac{1}{T} \int_{0}^{T} |\sum_{n \leq x} n^{it}|^{2q} dt \ll \Biggl(\frac{x}{1 + (1-q)\sqrt{\log\log(10L_{T})}} \Biggr)^q ,
\end{equation}
where\footnote{The appearances of $L_r$ and $L_T$ everywhere reflect a well known ``symmetry'' between $\sum_{n \leq x} \chi(n)$ and $\sum_{n \leq r/x} \chi(n)$, and between $\sum_{n \leq x} n^{it}$ and $\sum_{n \leq |t|/x} n^{it}$. This is explained in some detail in the author's paper~\cite{harpertypicalchar}, in Szab\'o's papers~\cite{szaboupper,szabolower}, and in the survey~\cite{harperbtsqc}. In the present paper, the size of $x$ with which we work means that we shall simply have $L_r, L_T = x$, and the ``symmetry'' will be irrelevant.} $L_T := \min\{x, T/x\}$. See the author's survey paper~\cite{harperbtsqc} for some general background and discussion of these issues.

Our goal here is to make progress on the lower bound problem for low moments. Again, on first inspection the shape of the bounds in \eqref{eqbtscupper} and \eqref{eqbtsccontupper} might seem very peculiar, and thus likely to be improvable. But motivated by order of magnitude results of the author~\cite{harperrmflow} in the model setting of random multiplicative functions, it was conjectured in \cite{harpertypicalchar} that \eqref{eqbtscupper} is sharp for all $x \leq 0.99r$, say. (Note that by periodicity, we have $|\sum_{n \leq x} \chi(n)| = |\sum_{n < r-x} \chi(n)|$ for non-principal characters $\chi$ mod $r$, allowing the situation where $x > 0.99r$ to be understood as well.) It was also conjectured~\cite{harpertypicalchar} that \eqref{eqbtsccontupper} is sharp for all $1 \leq x \leq T$. As discussed in \cite{harpertypicalchar}, the methods of that paper may be applicable to obtain sharp lower bounds when $x \leq e^{\log^{c}r}$ and $x \leq e^{\log^{c}T}$, respectively (for $c$ some small constant), but such an argument has not been worked out explicitly. The best existing lower bounds, available on the full range of $x$ for unweighted character sums $\frac{1}{r-1} \sum_{\chi \; \text{mod} \; r} |\sum_{n \leq x} \chi(n)|^{2q}$, and on the ranges $x \leq \sqrt{r}$ and $x \leq \sqrt{T}$ for the moments $\frac{1}{r-1} \sum_{\chi \; \text{mod} \; r} |\sum_{n \leq x} h(n) \chi(n)|^{2q}$ and $\frac{1}{T} \int_{0}^{T} |\sum_{n \leq x} n^{it}|^{2q} dt$, are due to La Bret\`eche, Munsch and Tenenbaum~\cite{bretechemunschten}. These differ from our upper bounds by powers of $\log L$ and $\log(L_T)$.

Our main result is:
\begin{thmlb1}\label{thmlb1}
Let $r$ be a large prime. Then uniformly for any $1 \leq x \leq r^{0.499}$ and any $0 \leq q \leq 1$, we have
$$ \frac{1}{r-1} \sum_{\chi \; \text{mod} \; r} |\sum_{n \leq x} \chi(n)|^{2q} \gg \Biggl(\frac{x}{1 + (1-q)\sqrt{\log\log(10x)}} \Biggr)^q . $$
\end{thmlb1}

In view of \eqref{eqbtscupper}, this is sharp, giving us order of magnitude results for these moments provided $1 \leq x \leq r^{0.499}$. There is no special significance to the exponent 0.499--- the proof will allow any {\em fixed} exponent $< 1/2$. More generally, we get the following sharp lower bounds:

\begin{thmlb2}\label{thmlb2}
Let $r$ be a large prime. Then uniformly for any $1 \leq x \leq r^{0.499}$, any $0 \leq q \leq 1$, and any multiplicative function $h(n)$ that has absolute value 1 on primes and absolute value at most 1 on prime powers, we have
$$ \frac{1}{r-1} \sum_{\chi \; \text{mod} \; r} |\sum_{n \leq x} h(n) \chi(n)|^{2q} \gg \Biggl(\frac{x}{1 + (1-q)\sqrt{\log\log(10x)}} \Biggr)^q . $$
\end{thmlb2}

\begin{thmlb3}\label{thmlb3}
Let $T$ be a large real number. Then uniformly for any $1 \leq x \leq T^{0.499}$ and any $0 \leq q \leq 1$, we have
$$ \frac{1}{T} \int_{0}^{T} |\sum_{n \leq x} n^{it}|^{2q} dt \gg \Biggl(\frac{x}{1 + (1-q)\sqrt{\log\log(10x)}} \Biggr)^q . $$
\end{thmlb3}

As will be seen, a restriction of the shape $x \leq r^{1/2 - o(1)}$ or $x \leq T^{1/2 - o(1)}$ arises naturally in our proofs. Nevertheless, as discussed above one expects that Theorems \ref{thmlb1} and \ref{thmlb3} ought to hold, with $\sqrt{\log\log(10x)}$ replaced by $\sqrt{\log\log(10L)}$ or $\sqrt{\log\log(10L_{T})}$, for all $x \leq 0.99r$ or $x \leq T$. We emphasise that for the moments $\frac{1}{r-1} \sum_{\chi \; \text{mod} \; r} |\sum_{n \leq x} h(n) \chi(n)|^{2q}$ involving a multiplicative twist $h(n)$, we do {\em not} generally expect the upper bound we have (as in \eqref{eqbtscupper}) to be sharp for $x$ close to $r$, since in general there is no reason why $\log\log(10L)$ should appear in place of the stronger saving $\log\log(10x)$. (But in specific cases, such as when $h(n) \equiv 1$, ``symmetry'' means that the appearance of $\log\log(10L)$ is indeed required.) For example, in the important case where $h$ is the M\"obius function $\mu(n)$, the author~\cite{harpertypicalchar} conjectured that we should have $\frac{1}{r-1} \sum_{\chi \; \text{mod} \; r} |\sum_{n \leq x} \mu(n) \chi(n)|^{2q} \ll_{A} (\frac{x}{1 + (1-q)\sqrt{\log\log(10x)}} )^q$ for all $x \leq r^A$ and {\em any} fixed $A > 0$, and this has been proved\footnote{Strictly speaking, Wang and Xu~\cite{wangxu} study the moments $\frac{1}{r-1} \sum_{\chi \; \text{mod} \; r} |\sum_{n \leq x} \lambda(n) \chi(n)|^{2q}$ twisted by the Liouville function $\lambda(n)$ rather than the M\"obius function.} by Wang and Xu~\cite{wangxu} assuming certain other strong conjectures.

In a companion paper~\cite{harperlowchar2}, we extend our proofs (and exploit the special properties of unweighted character sums) to establish an analogue of Theorem \ref{thmlb1}, with $\sqrt{\log\log(10x)}$ replaced by $\sqrt{\log\log(10L)}$, for all $x \leq 0.99r$. This gives us a full order of magnitude understanding of the unweighted moments $\frac{1}{r-1} \sum_{\chi \; \text{mod} \; r} |\sum_{n \leq x} \chi(n)|^{2q}$ for $q \leq 1$ and all $x \leq 0.99r$ (and thus for all $x$). In particular, as described in \cite{harpertypicalchar,harperbtsqc}, these types of estimates with $x \approx \sqrt{r}$ may be applied to obtain new non-vanishing results for Dirichlet theta functions $\theta(1,\chi)$ (since these are essentially mildly weighted character sums of length $\approx \sqrt{r}$).

\subsection{Ideas from the proofs}\label{subsecproofideas}
The majority of this paper is devoted to the proof of Theorem \ref{thmlb1}. The proofs of Theorems \ref{thmlb2} and \ref{thmlb3} are fairly straightforward variants of this, and will be discussed briefly at the end.

Throughout we shall let $r$ denote a large prime modulus, and for concision we shall write $\Echar$ to denote averaging over all Dirichlet characters modulo $r$. Thus if $W(\chi)$ is any function, then $\Echar W := \frac{1}{r-1} \sum_{\chi \; \text{mod} \; r} W(\chi)$. We also recall that a Steinhaus random multiplicative function $f(n)$ is obtained by letting $(f(p))_{p \; \text{prime}}$ be a sequence of independent Steinhaus random variables (i.e. distributed uniformly on the unit circle $\{|z|=1\}$), and then setting $f(n) := \prod_{p^{a} || n} f(p)^{a}$ for all $n \in \N$, where $p^a || n$ means that $p^a$ is the highest power of the prime $p$ that divides $n$. These satisfy the orthogonality relation $\E f(n) \overline{f(m)} = \textbf{1}_{n=m}$ for all $n,m \in \N$, so in particular (using the primality of $r$ a little) we get
\begin{equation}\label{eqnorthcharrmf}
\Echar \chi(n) \overline{\chi(m)} = \textbf{1}_{n=m} = \E f(n) \overline{f(m)} \;\;\;\;\;\;\;\; \forall \; 1 \leq n, m < r .
\end{equation}

In this subsection, we shall try to outline the main ideas in the proof of Theorem \ref{thmlb1}, and reduce this to proving two other Key Propositions. In the next subsection, we make some further remarks including connections with the wider literature.

\vspace{12pt}
Our basic strategy is a staple of analytic number theory, being the same approach taken in most modern work on (lower bounds for) {\em moments of $L$-functions}, as well as the methods of {\em amplification}, {\em mollification}, and {\em resonance}. That is, one tries to define a proxy object that mimics the relevant property of the object of study, but is easier to analyse. Then one computes some form of correlation of the proxy object with the main object, as well as an average bound for the size of the proxy, and compares the two. For example, when studying moments of $L$-functions one chooses a proxy that mimics some relevant fixed power of the $L$-function; in amplification and resonance, the proxy boosts the contribution from large values, for example by (somewhat) mimicking a high power of the main object; and in mollification, the proxy (somewhat) mimics the inverse of the main object, so that the product of the two should have roughly constant size on average. The key issue then becomes the choice of the proxy object. In the settings mentioned, the proxy is generally chosen as a short Dirichlet polynomial with appropriate {\em combinatorial coefficients}, for example something like generalised divisor functions in the moment and resonance situations, and something like the (smoothed) M\"obius function in the mollification setting. Here, our choice of proxy $I(\chi)$ will have a somewhat different shape.

\vspace{12pt}
As motivation, notice Perron's formula implies $\sum_{n \leq x} \chi(n) \approx \frac{1}{2\pi i} \int_{1/2 - i\infty}^{1/2 + i\infty} (\sum_{n \leq x} \frac{\chi(n)}{n^s}) \frac{x^s}{s} ds$. Taking the line of integration to have real part 1/2 seems sensible, since we are looking to show that the usual size of $\sum_{n \leq x} \chi(n)$ is around (indeed, slightly smaller than) $\sqrt{x}$. We shall produce $I(\chi)$ by editing this Perron integral in three ways. The first edit, which is standard and hopefully requires little further explanation, is that we shall replace the infinite integral by a suitable finite truncation (up to height tending slowly to infinity with $x$). The second edit is that we shall replace $\sum_{n \leq x} \frac{\chi(n)}{n^s}$ by a related, but more nicely structured, Dirichlet polynomial. We will discuss this further in a moment. The third (and most important) edit is that we shall insert a ``barrier'' into the integral, discarding (or greatly penalising) points $s$ where the behaviour of $\sum_{n \leq x} \frac{\chi(n)}{n^s}$ (or the more structured sum we actually use) is ``abnormal''. The choice of barrier is characteristic of the {\em multiplicative chaos} phenomena underlying bounds like \eqref{eqbtscupper} and \eqref{eqbtsccontupper}, and will be crucial for our calculations to succeed. Again, we shall discuss this point in more detail below.

Let us now define
\begin{eqnarray}
I(\chi) = I_{x,q}(\chi) & := & \frac{1}{2\pi i} \int_{1/2-i\mathcal{T}}^{1/2+i\mathcal{T}} \Biggl(\sum_{x^{1-3\beta} < p \leq x^{1-\beta}} \frac{\chi(p)}{p^s} \Biggr) F_{\chi}^{*}(s) G_{\chi}(s) \frac{x^s}{s} ds \nonumber \\
& = & \frac{\sqrt{x}}{2\pi} \int_{-\mathcal{T}}^{\mathcal{T}} \Biggl(\sum_{x^{1-3\beta} < p \leq x^{1-\beta}} \frac{\chi(p)}{p^{1/2+iv}} \Biggr) F_{\chi}^{*}(1/2 + iv) G_{\chi}(1/2 + iv) \frac{x^{iv}}{1/2 + iv} dv , \nonumber
\end{eqnarray}
where $p$ runs over primes. Here we let
$$ F_{\chi}^{*}(s) := \Biggl(\sum_{\substack{n \leq x^{\beta}, \\ n \; \text{is} \; P \; \text{smooth}}} \frac{\chi(n)}{n^s} \Biggr) \cdot \Biggl(\sum_{\substack{n \leq x^{2\beta}, \\ n \; \text{is} \; P \; \text{rough}}} \frac{\chi(n)}{n^s} \Biggr) , $$
recalling that a number $n$ is said to be $P$-smooth if all of its prime factors are $\leq P$, and to be $P$-rough if all of its prime factors are $> P$. And $G_{\chi}(s) = G_{\chi, q,P}(s)$ will be (a Dirichlet polynomial approximation to a smoothed version of) a barrier, placing simultaneous restrictions on the sizes of the partial sums $\Re \sum_{p \leq P^{e^{-j}}} \sum_{k=1,2} \frac{\chi(p^k)}{kp^{ks}}$ for various $j$, and thus on the sizes of the partial Euler products $\prod_{p \leq P^{e^{-j}}} |1 - \frac{\chi(p)}{p^{s}}|^{-1}$. See section \ref{secpasstormf}, below, for our precise choice of $G_{\chi}(s)$. In our notation $I(\chi)$ we suppress mention of the parameters $\beta$ (a small positive quantity), and $P,\mathcal{T}$ (which will have size $x^{o(1)}$), on which $I(\chi)$ of course also depends. Note that $I(\chi)$ only depends on $q$ through the choice of barrier $G_{\chi, q,P}(s)$, as $q$ approaches 1 (and so we get closer to the trivial problem of calculating $\frac{1}{r-1} \sum_{\chi \; \text{mod} \; r} |\sum_{n \leq x} \chi(n)|^{2}$, and the saving factor $1 + (1-q)\sqrt{\log\log x}$ shrinks) the barrier becomes more relaxed.

Theorem \ref{thmlb1} is an easy consequence of the following correlation and moment bounds.

\begin{keyprop1}\label{keyprop1}
Uniformly for all large $x$, all $0 \leq q \leq 1$, and all small $\beta \geq 1/\log^{0.1}x$ (say) such that $x^{1+2\beta+ \beta/(\log\log x)^{10}} < r$, the following is true. For $I(\chi)$ defined as above, and with $\mathcal{T},P$ chosen as in \eqref{eqselectTP} below, we have
$$ \Echar |\sum_{n \leq x} \chi(n)| |I(\chi)| \gg \frac{\beta x}{1 + (1-q)\sqrt{\log\log x}} . $$
\end{keyprop1}

\begin{keyprop2}\label{keyprop2}
Uniformly for all large $x$, all $0 \leq q \leq 1$, and all small $\beta \geq 1/\log^{0.1}x$ (say) such that $x^{2+4\beta+ \beta/(\log\log x)^{10}} < r$, the following is true. For $I(\chi)$ defined as above, and with $\mathcal{T},P$ chosen as in \eqref{eqselectTP} below, we have
$$ \Echar |I(\chi)|^2 \ll \frac{\beta x}{1 + (1-q)\sqrt{\log\log x}} , \;\;\; \Echar |I(\chi)|^4 \ll e^{2\min\{\sqrt{\log\log x}, \frac{1}{1-q}\}} \left( \frac{\beta x}{1 + (1-q)\sqrt{\log\log x}} \right)^2 . $$
\end{keyprop2}

\begin{proof}[Proof of Theorem \ref{thmlb1}, assuming Key Propositions \ref{keyprop1} and \ref{keyprop2}]
We may assume that $x$ is large, otherwise Theorem \ref{thmlb1} is trivial since e.g. we always have the bound $\Echar |\sum_{n \leq x} \chi(n)|^{2q} \geq x^{-2(1-q)} \Echar |\sum_{n \leq x} \chi(n)|^{2} = x^{-2(1-q)} \lfloor x \rfloor$, and this is $\gg x^q$ for bounded $x$.

Using H\"older's inequality, for any $2/3 \leq q \leq 1$ we get
\begin{eqnarray}
&& \Echar |\sum_{n \leq x} \chi(n)| |I(\chi)| = \Echar |\sum_{n \leq x} \chi(n)| |I(\chi)|^{(3q-2)/q} |I(\chi)|^{2(1-q)/q} \nonumber \\
& \leq & \left( \Echar |\sum_{n \leq x} \chi(n)|^{2q} \right)^{1/2q}  \left( \Echar |I(\chi)|^{2} \right)^{(3q-2)/2q} \left( \Echar |I(\chi)|^{4} \right)^{(1-q)/2q} . \nonumber
\end{eqnarray}
Fix $\beta$ as a small constant such that the conditions $x^{1+2\beta+ \beta/(\log\log x)^{10}}, x^{2+4\beta+ \beta/(\log\log x)^{10}} < r$ are both satisfied (which crucially is possible under our assumption that $x \leq r^{0.499}$). Using Key Proposition \ref{keyprop1} to lower bound the left hand side, and Key Proposition \ref{keyprop2} to upper bound the latter two terms on the right, we deduce
\begin{eqnarray}
\Echar |\sum_{n \leq x} \chi(n)|^{2q} & \geq & \frac{(\Echar |\sum_{n \leq x} \chi(n)| |I(\chi)|)^{2q}}{\left( \Echar |I(\chi)|^{2} \right)^{3q-2} \left( \Echar |I(\chi)|^{4} \right)^{1-q}} \nonumber \\
& \gg & \frac{(\frac{\beta x}{1 + (1-q)\sqrt{\log\log x}})^{2q}}{(\frac{\beta x}{1 + (1-q)\sqrt{\log\log x}})^{3q-2} (e^{2\min\{\sqrt{\log\log x}, \frac{1}{1-q}\}} (\frac{\beta x}{1 + (1-q)\sqrt{\log\log x}})^2 )^{1-q} } \nonumber \\
& \gg & (\frac{\beta x}{1 + (1-q)\sqrt{\log\log x}})^{q} \gg (\frac{x}{1 + (1-q)\sqrt{\log\log x}})^{q} . \nonumber
\end{eqnarray}
This proves the theorem for all $2/3 \leq q \leq 1$.

For smaller $q$, we can deduce our desired result with a further easy application of H\"older's inequality. For all $0 \leq q \leq 2/3$, we have
\begin{eqnarray}
\Echar |\sum_{n \leq x} \chi(n)|^{4/3} & = & \Echar |\sum_{n \leq x} \chi(n)|^{q/(3(3/2 - 2q))} |\sum_{n \leq x} \chi(n)|^{(2-3q)/(3/2 - 2q)} \nonumber \\
& \leq & \left( \Echar |\sum_{n \leq x} \chi(n)|^{2q} \right)^{1/(6(3/2 - 2q))} \left( \Echar |\sum_{n \leq x} \chi(n)|^{3/2} \right)^{(4-6q)/(3(3/2 - 2q))} . \nonumber
\end{eqnarray}
Using the fact that $\Echar |\sum_{n \leq x} \chi(n)|^{4/3} \gg \left( \frac{x}{\sqrt{\log\log x}} \right)^{2/3}$, along with the upper bound $\Echar |\sum_{n \leq x} \chi(n)|^{3/2} \ll \left( \frac{x}{\sqrt{\log\log x}} \right)^{3/4}$ from \eqref{eqbtscupper}, we then find
$$ \Echar |\sum_{n \leq x} \chi(n)|^{2q} \geq \frac{(\Echar |\sum_{n \leq x} \chi(n)|^{4/3} )^{3(3-4q)}}{( \Echar |\sum_{n \leq x} \chi(n)|^{3/2} )^{4(2-3q)}} \gg \left( \frac{x}{\sqrt{\log\log x}} \right)^{q} . $$
\end{proof}

It seems worth observing that we only need to deploy the upper bound for $\Echar |I(\chi)|^2$ because we prove Theorem \ref{thmlb1} uniformly for $q$ approaching 1. For $q$ bounded away from 1, a simpler application of H\"older's inequality could be made, only relying on the lower bound for $\Echar |\sum_{n \leq x} \chi(n)| |I(\chi)|$ and the upper bound for $\Echar |I(\chi)|^4$.

We also note the behaviour of the various quantities in Key Propositions \ref{keyprop1} and \ref{keyprop2} as $q$ approaches 1, and our barrier (encoded in $G_{\chi, q,P}(s)$) relaxes. Both $\Echar |\sum_{n \leq x} \chi(n)| |I(\chi)|$ and $\Echar |I(\chi)|^2$ grow in the way one might expect given the statement of Theorem \ref{thmlb1}--- note that our estimates for both are of the same order, providing some support that $I(\chi)$ is indeed a reasonable proxy for $\sum_{n \leq x} \chi(n)$. The fourth moment $\Echar |I(\chi)|^4$ grows like the square of the second moment (as one would hope), multiplied by the factor $e^{2\min\{\sqrt{\log\log x}, \frac{1}{1-q}\}}$ which is negligible when raised to the power $(1-q)/2q$. Our choice of $G_{\chi, q,P}(s)$ is calibrated to control this blow-up factor produced in the fourth moment estimate (see the discussion below)--- we make the barrier as relaxed as possible subject to the blow-up remaining negligible when we apply H\"older's inequality (as above), and the optimally relaxed choice ultimately produces the $(1-q)$ multiplier in all of the denominators (and our main theorems).

\vspace{12pt}
We shall now discuss the form of $I(\chi)$ somewhat further.

Compared with a naive Perron integral, we have replaced $\sum_{n \leq x} \frac{\chi(n)}{n^s}$ with the product $(\sum_{x^{1-3\beta} < p \leq x^{1-\beta}} \frac{\chi(p)}{p^s} ) F_{\chi}^{*}(s) = (\sum_{x^{1-3\beta} < p \leq x^{1-\beta}} \frac{\chi(p)}{p^s} ) \cdot (\sum_{\substack{n \leq x^{\beta}, \\ n \; \text{is} \; P \; \text{smooth}}} \frac{\chi(n)}{n^s} ) \cdot (\sum_{\substack{n \leq x^{2\beta}, \\ n \; \text{is} \; P \; \text{rough}}} \frac{\chi(n)}{n^s} )$. Beyond the specific shape of this, two things to note are the product form, and crucially the lengths of the involved Dirichlet polynomials.

The summands in the three Dirichlet polynomials only depend (thanks to total multiplicativity of Dirichlet characters) on the values of $\chi$ on disjoint sets of primes, namely those on the ranges $(x^{1-3\beta}, x^{1-\beta}], [2,P], (P, x^{2\beta}]$. Together with the fact that there is no interaction between the ranges in the three summations, this essentially means that we will be able to analyse these Dirichlet polynomials one at a time in our computations. (On the random multiplicative side, this corresponds to the fact that we can condition on the behaviour of $f$ on some primes without affecting its behaviour on any others.)

The total length of our sums here is $x^{1-\beta} \cdot x^{\beta} \cdot x^{2\beta} = x^{1+2\beta}$. The ``barrier'' $G_{\chi, q,P}(s)$ is a short Dirichlet polynomial of length $x^{o(1)}$, more specifically this turns out to have length roughly $x^{\beta/(\log\log x)^{10}}$. This means that when computing the second moment type quantities $\Echar |\sum_{n \leq x} \chi(n)| |I(\chi)|, \Echar |I(\chi)|^2$, by \eqref{eqnorthcharrmf} the character average will match the corresponding random multiplicative function average (or, equivalently, only diagonal terms are relevant) provided $x^{1+2\beta+ \beta/(\log\log x)^{10}} < r$. For the fourth moment $\Echar |I(\chi)|^4$, the lengths are squared when everything is expanded out and we require $x^{2+4\beta+ \beta/(\log\log x)^{10}} < r$ to pass to the random multiplicative side. Notice these are the conditions of Key Propositions \ref{keyprop1} and \ref{keyprop2}.

These requirements, in particular the restriction $x^{2+4\beta+ \beta/(\log\log x)^{10}} < r$ coming from the fourth moment of our proxy $I(\chi)$, restrict us to working with $x$ a little smaller than $\sqrt{r}$. The lengths of the Dirichlet polynomials defining $I(\chi)$ could be adjusted, but it seems unreasonable to expect a proxy for $\sum_{n \leq x} \chi(n)$ having length much less than $x$, especially since we analyse $\Echar |\sum_{n \leq x} \chi(n)| |I(\chi)|$ by estimating $\Echar \sum_{n \leq x} \chi(n) \overline{I(\chi)}$ (without absolute values) so some genuine reinforcement between $\sum_{n \leq x} \chi(n)$ and $I(\chi)$ is required. The crucial barrier necessarily adds something to the lengths of the Dirichlet polynomials, so constructing $I(\chi)$ with length slightly larger than $x$ (as we do) seems more or less the best that we can hope for. The Dirichlet polynomial lengths in the fourth moment type calculation is the key issue that must be addressed and worked around in the companion paper~\cite{harperlowchar2}, to extend the range of $x$ in Theorem \ref{thmlb1}.

\vspace{12pt}
Each Dirichlet polynomial $\sum_{x^{1-3\beta} < p \leq x^{1-\beta}} \frac{\chi(p)}{p^s} , \sum_{\substack{n \leq x^{\beta}, \\ P \; \text{smooth}}} \frac{\chi(n)}{n^s} , \sum_{\substack{n \leq x^{2\beta}, \\ P \; \text{rough}}} \frac{\chi(n)}{n^s}$ plays a different and important role in our calculations. They could be modified, but this particular combination seems to work quite nicely. This is perhaps easiest to think about on the random multiplicative function side, although the calculations are (essentially) the same. The Dirichlet polynomial over large primes, $\sum_{x^{1-3\beta} < p \leq x^{1-\beta}} \frac{f(p)}{p^s}$ or $\sum_{x^{1-3\beta} < p \leq x^{1-\beta}} \frac{\chi(p)}{p^s}$, gives us summands $f(p)$ that behave independently of one another and have easily understood interactions with everything else (since these primes only appear once inside our Dirichlet polynomials, and at most one of them divides any $n$ in the sum $\sum_{n \leq x} f(n)$ or $\sum_{n \leq x} \chi(n)$). This ultimately means that the calculation of the second moment quantities $\E \sum_{n \leq x} f(n) \overline{I(f)}, \E |I(f)|^2$ is reduced to handling $\frac{x}{\log x} \E \int_{-\mathcal{T}}^{\mathcal{T}} \frac{|F_{f}^{*}(\frac{1}{2} + iv)|^2 |G_{f}(\frac{1}{2} + iv)|^2}{|1/2 + iv|^2} dv$, (notice this is obviously real and non-negative, which is certainly not obvious a priori for $\Echar \sum_{n \leq x} \chi(n) \overline{I(\chi)}$ or $\E \sum_{n \leq x} f(n) \overline{I(f)}$), and calculation of the fourth moment $\E |I(f)|^4$ reduces to handling $\E(\frac{x}{\log x} \int_{-\mathcal{T}}^{\mathcal{T}} \frac{|F_{f}^{*}(\frac{1}{2} + iv)|^2 |G_{f}(\frac{1}{2} + iv)|^2}{|1/2 + iv|^2} dv )^2$. See especially the application of Khintchine's inequality at the beginning of section \ref{subsecfourthkh}, below, as well as the calculations in section \ref{subsecinitcond}.

With an appropriate choice of the parameter $P$ relative to $x^{\beta}$, the smooth Dirichlet polynomial $\sum_{\substack{n \leq x^{\beta}, \\ P \; \text{smooth}}} \frac{f(n)}{n^s}$ may be replaced, with acceptable average error, by the Euler product $F_{P}(s) = \prod_{p \leq P} (1 - \frac{f(p)}{p^s})^{-1}$. Note that we could not directly include an Euler product in the definition of $I(\chi)$, since when expanded out this would lead to Dirichlet polynomials far too long to invoke \eqref{eqnorthcharrmf}. The Euler product (a product of independent terms), together with the barrier $G_{\chi, q,P}(s)$ or $G_{f, q,P}(s)$ constraining its growth (see below), is the ultimate source of the factor $1 + (1-q)\sqrt{\log\log(10x)}$ in our bounds.

The Dirichlet polynomial over $P$-rough numbers plays the least crucial role in our analysis, but several comments about it are still in order. To possibly obtain the sharp lower bound in Theorem \ref{thmlb1}, the Dirichlet polynomials we use to define $I(\chi)$ must, between them, sum over a positive density subset of integers (thereby correlating, potentially, with a positive density subsum of $\sum_{n \leq x} \chi(n)$). Those integers with a large prime factor $x^{1-3\beta} < p \leq x^{1-\beta}$, and which are otherwise $P$-smooth, would {\em not} be a positive density subset--- allowing a $P$-rough component as well rectifies this. Indeed, the factor $\beta$ in our lower bound in Key Proposition \ref{keyprop1} directly corresponds to the density of the set of integers we are summing over, see the calculations at the end of section \ref{subsecreducesingle}. Several of our calculations involve extracting cancellation from appropriate sums over $P$-rough numbers, the philosophy being that terms $p^{-ih}$ with $p > P$ oscillate rapidly. In places this is combined with additional cancellation coming from sums over primes $x^{1-3\beta} < p \leq x^{1-\beta}$, here we benefit from the fact that taking a product of Dirichlet polynomials (multiplicative convolution) multiplies the savings from the different components. Finally, we emphasise that although $P$ must be rather smaller than $x^{2\beta}$, it is also important that it is {\em not too much smaller}. In some of the calculations in section \ref{subsecexploitbarr}, we incur a loss from the $P$-rough sums because the barrier $G_{\chi, q,P}(s)$ does not constrain their growth (it cannot, otherwise it could not be constructed as a sufficiently short Dirichlet polynomial). We have some saving coming from the barrier on the $P$-smooth sums to compensate (indeed this is the whole point of the barrier, see below), but this is only adequate if the loss from the $P$-rough sums, which is of the shape $(\sum_{\substack{n \leq x^{2\beta}, \\ P \; \text{rough}}} \frac{1}{n})^{O(1)} \asymp (\frac{\log(x^{2\beta})}{\log P})^{O(1)}$, is not too large.

\vspace{12pt}
To explain the vital role of the barrier\footnote{It would be more accurate to describe $G_{\chi}(s)$ as a Dirichlet polynomial approximation to a barrier. Even more precisely, $G_{\chi}(s)$ is an approximation to a smoothed version $W_{\chi}(s)$ of a barrier, which we ultimately analyse by (up to acceptable errors) sandwiching it between two ``sharp'' barriers $L_{\chi}$ and $U_{\chi}$. We suppress all of these subtleties for our introductory discussion, see section \ref{secpasstormf} for the details.} $G_{\chi}(s) = G_{\chi, q,P}(s)$, consider the problems of bounding $\Echar |I(\chi)|^4$ or $\E |I(f)|^4$ (which, by \eqref{eqnorthcharrmf} and our condition $x^{2+4\beta+ \beta/(\log\log x)^{10}} < r$, are in fact the same problem). As noted earlier, exploiting the presence of the Dirichlet polynomial $\sum_{x^{1-3\beta} < p \leq x^{1-\beta}} \frac{f(p)}{p^s}$ over large primes, followed by a simple smoothing argument, one can show fairly readily that
\begin{equation}
\E|I(f)|^4 \ll x^2 \E\left( \frac{1}{\log x} \int_{-\mathcal{T}}^{\mathcal{T}} \frac{|F_{f}^{*}(\frac{1}{2} + iv)|^2 |G_{f}(\frac{1}{2} + iv)|^2}{|1/2 + iv|^2} dv \right)^2 . \nonumber
\end{equation}
See display \eqref{eqfourthafterred}, below. Thanks to the denominator $|1/2 + iv|^2 = 1/4 + v^2$, the length of the integral is fairly immaterial here and it essentially suffices to handle the portion where $|v| \leq 1/2$, say. As also noted earlier, with appropriate parameter choices the Dirichlet polynomial $\sum_{\substack{n \leq x^{\beta}, \\ P \; \text{smooth}}} \frac{f(n)}{n^s}$ inside $F_{f}^{*}(s)$ may be replaced by the random Euler product $F_{P}(s)$. So we are left to understand something like
$$ \E\Biggl( \frac{1}{\log x} \int_{-1/2}^{1/2} |\sum_{\substack{n \leq x^{2\beta}, \\ n \; \text{is} \; P \; \text{rough}}} \frac{f(n)}{n^{1/2 + iv}}|^2 |F_{P}(\frac{1}{2} + iv)|^2 |G_{f}(\frac{1}{2} + iv)|^2 dv \Biggr)^2 . $$
Expanding the square, this requires controlling correlations $\E |\sum_{\substack{n \leq x^{2\beta}, \\ n \; \text{is} \; P \; \text{rough}}} \frac{f(n)}{n^{1/2 + iv}}|^2 |F_{P}(\frac{1}{2} + iv)|^2 |G_{f}(\frac{1}{2} + iv)|^2 |\sum_{\substack{n \leq x^{2\beta}, \\ n \; \text{is} \; P \; \text{rough}}} \frac{f(n)}{n^{1/2 + iw}}|^2 |F_{P}(\frac{1}{2} + iw)|^2 |G_{f}(\frac{1}{2} + iw)|^2$ for all $|v|, |w| \leq 1/2$. Note that the contributions from $P$-smooth and $P$-rough numbers are independent.

In the most optimistic circumstance, we could hope that the $P$-rough sums at $v$ and $w$ would decorrelate (since terms $p^{-i(v-w)}$ with $p > P$ oscillate rapidly), thus contributing $\approx  (\E |\sum_{\substack{n \leq x^{2\beta}, \\ n \; \text{is} \; P \; \text{rough}}} \frac{f(n)}{n^{1/2 + iv}}|^2) (\E |\sum_{\substack{n \leq x^{2\beta}, \\ n \; \text{is} \; P \; \text{rough}}} \frac{f(n)}{n^{1/2 + iw}}|^2) = (\sum_{\substack{n \leq x^{2\beta}, \\ n \; \text{is} \; P \; \text{rough}}} \frac{1}{n})^2 \asymp (\frac{\log(x^{2\beta})}{\log P})^{2}$. To establish the fourth moment bound claimed in Key Proposition \ref{keyprop2}, we would then like to show that on average over $v,w$, the $P$-smooth contribution is (a little smaller than) $\log^{2}P$. But without any barriers $G_{f}(\frac{1}{2} + iv), G_{f}(\frac{1}{2} + iw)$, it is easy to calculate that
\begin{equation}\label{eqnnobarrcorr}
\E |F_{P}(\frac{1}{2} + iv)|^2 |F_{P}(\frac{1}{2} + iw)|^2 \asymp (\log P)^{2} \min\{(\log P)^2 , \frac{1}{|v-w|^2 }\} \;\;\;\;\; \forall \; |v|, |w| \leq 1/2 .
\end{equation}
Integrating this over $v,w$, we would not obtain anything close to $\log^{2}P$, but the unacceptably large $\log^{3}P$. 

{\em The purpose of the barrier is to suppress this blow-up created by ``nearby'' $v,w$.} The phenomena at work here, which we shall try to succinctly explain, underlie all bounds involving so-called (critical) multiplicative chaos behaviour. See e.g. section 3 of the author's survey paper~\cite{harperrmf3} for much more discussion of these issues, in the context of random multiplicative functions.

In the computation of $\E |F_{P}(\frac{1}{2} + iv)|^2 |F_{P}(\frac{1}{2} + iw)|^2$, one finds that the parts of the Euler products on primes $p \leq e^{1/|v-w|}$ (or $p \leq P$, if $|v-w| \leq 1/\log P$) are highly correlated, because for such $p$ there is little oscillation in the factors $p^{-i(v-w)} = e^{-i(v-w)\log p}$ that could cause the products at $v$ and $w$ to behave differently. On the other hand, on $e^{1/|v-w|} < p \leq P$ there is sufficient oscillation that one does have the desired decorrelation $\E \prod_{e^{1/|v-w|} < p \leq P} |1 - \frac{f(p)}{p^{1/2 + iv}}|^{-2} \prod_{e^{1/|v-w|} < p \leq P} |1 - \frac{f(p)}{p^{1/2 + iw}}|^{-2} \asymp (\E \prod_{e^{1/|v-w|} < p \leq P} |1 - \frac{f(p)}{p^{1/2 + iv}}|^{-2}) (\E \prod_{e^{1/|v-w|} < p \leq P} |1 - \frac{f(p)}{p^{1/2 + iw}}|^{-2})$. So we want to set things up so that our barriers $G_{f}(\frac{1}{2} + iv), G_{f}(\frac{1}{2} + iw)$ appropriately restrict the sizes of the partial Euler products on $p \leq e^{1/|v-w|}$. Since the specific size of $|v-w|$ (varying with both $v$ and $w$) cannot be built into $G_{f}(\frac{1}{2} + iv)$ or $G_{f}(\frac{1}{2} + iw)$, this leads us to construct $G_{f}(\frac{1}{2} + iv)$ that simultaneously restricts {\em all} partial products $\prod_{p \leq P^{e^{-j}}} |1 - \frac{f(p)}{p^{1/2 + iv}}|^{-1}$. The terminology ``barrier'' may be explained on viewing the sequence of partial products, for varying $j$, as a (multiplicative) walk restricted to lie within certain bounds.

In order to control the blow-up, it turns out that requiring $\prod_{p \leq P^{e^{-j}}} |1 - \frac{f(p)}{p^{1/2 + iv}}|^{-1}$ to be bounded by $\log(P^{e^{-j}})$ is roughly sufficient. For then, invoking this barrier with $e^j \approx 1 + |v-w|\log P$ to bound the subproduct of $F_{P}(\frac{1}{2} + iv)$ up to $\approx \min\{e^{1/|v-w|}, P\}$, we can bound $\E |F_{P}(\frac{1}{2} + iv)|^2 |G_{f}(\frac{1}{2} + iv)|^2 |F_{P}(\frac{1}{2} + iw)|^2 |G_{f}(\frac{1}{2} + iw)|^2$ by
$$ \approx \min\{(\log P)^2 , \frac{1}{|v-w|^2 }\} \E \prod_{e^{1/|v-w|} < p \leq P} |1 - \frac{f(p)}{p^{1/2 + iv}}|^{-2} |G_{f}(\frac{1}{2} + iv)|^2 |F_{P}(\frac{1}{2} + iw)|^2 |G_{f}(\frac{1}{2} + iw)|^2 . $$
Roughly speaking, the subproduct of $F_{P}(\frac{1}{2} + iw)$ up to $\min\{e^{1/|v-w|}, P\}$ is independent of everything remaining inside the expectation, so may be replaced by its mean square, which is $\asymp \min\{ \log P , \frac{1}{|v-w|}\}$. And the remaining products over large primes $e^{1/|v-w|} < p \leq P$ decorrelate, leaving us with
\begin{eqnarray}
& \approx & \min\{(\log P)^3 , \frac{1}{|v-w|^3 }\} (\E \prod_{e^{1/|v-w|} < p \leq P} |1 - \frac{f(p)}{p^{1/2 + iv}}|^{-2} |G_{f}(\frac{1}{2} + iv)|^2) \cdot \nonumber \\
&& \cdot (\E \prod_{e^{1/|v-w|} < p \leq P} |1 - \frac{f(p)}{p^{1/2 + iw}}|^{-2} |G_{f}(\frac{1}{2} + iw)|^2 ) . \nonumber
\end{eqnarray}
(This is only roughly true because $G_{f}(\frac{1}{2} + iv), G_{f}(\frac{1}{2} + iw)$ depend on the values of $f$ on all primes, including the small primes, which slightly interferes with the independence and decorrelation.) As a first approximation, if we now ignore the barriers $G_f$ we find $\E \prod_{e^{1/|v-w|} < p \leq P} |1 - \frac{f(p)}{p^{1/2 + iv}}|^{-2} = \E \prod_{e^{1/|v-w|} < p \leq P} |1 - \frac{f(p)}{p^{1/2 + iw}}|^{-2} \asymp \max\{1, |v-w|\log P\}$, giving an overall contribution $\approx (\log P)^{2} \min\{ \log P , \frac{1}{|v-w|}\}$. This may be compared with \eqref{eqnnobarrcorr}. Integrating over $|v|, |w| \leq 1/2$ would yield $(\log P)^2 \log\log P$, narrowly missing the size (slightly smaller than) $\log^{2}P$ that we are seeking.

To perfect this argument, one actually chooses the barrier to impose a slightly different upper bound on  $\prod_{p \leq P^{e^{-j}}} |1 - \frac{f(p)}{p^{1/2 + iv}}|^{-1}$ in place of $\log(P^{e^{-j}})$. It essentially suffices to have a bound $\frac{\log(P^{e^{-j}})}{(\log\log(P^{e^{-j}}))^{O(1)}}$, the doubly logarithmic factor killing off the $\log\log P$ divergence in our calculation above. Keeping the barriers $G_f$ in the final computation of $\E \prod_{e^{1/|v-w|} < p \leq P} |1 - \frac{f(p)}{p^{1/2 + iv}}|^{-2} |G_{f}(\frac{1}{2} + iv)|^2, \E \prod_{e^{1/|v-w|} < p \leq P} |1 - \frac{f(p)}{p^{1/2 + iw}}|^{-2} |G_{f}(\frac{1}{2} + iw)|^2$ is responsible for a further multiplier $\left( \frac{1}{1 + (1-q)\sqrt{\log\log P}} \right)^2$, which is $\asymp \left( \frac{1}{1 + (1-q)\sqrt{\log\log x}} \right)^2$ with our choice of the parameter $P$. See Probability Results \ref{probres4} and \ref{probres5}, below, and the calculations with them in sections \ref{subseclowercon} and \ref{subsecexploitbarr}. The presence of $\log\log P$ reflects the number of ``scales'' $P^{e^{-j}}$ on which the barrier constrains our multiplicative walk.

There are a couple of further refinements that we have already briefly discussed, and shall simply mention again. Firstly, our optimistic claim that the $P$-rough contributions (removed at the beginning of the above discussion) simply decorrelate turns out to be true when $|v-w| \geq 1/\log P$, but not when $v,w$ are closer together. See Probability Result \ref{probres3}, below. This generates a small additional blow-up, but with our parameter choices this can also be removed using the $(\log\log(P^{e^{-j}}))^{O(1)}$ denominator in our perfected barrier (note that when $|v-w| < 1/\log P$ we apply the barrier with $j=0$). Secondly, our actual choice of barrier also depends on the underlying moment exponent $q$, including a relaxation factor multiplying the bound $\frac{\log(P^{e^{-j}})}{(\log\log(P^{e^{-j}}))^{O(1)}}$ that increases as $q$ approaches 1. This is responsible for the (unhelpful, but by our choices harmless) multiplier $e^{2\min\{\sqrt{\log\log x}, \frac{1}{1-q}\}}$ in our upper bound for $\Echar |I(\chi)|^4$, but importantly also for the factor $(1-q)$ multiplying $\sqrt{\log\log x}$ in the denominators, sharpening our overall lower bounds when $q$ approaches 1.

\vspace{12pt}
In section \ref{secpasstormf}, we make our precise choice of the barrier $G_{\chi}(s)$, and perform some preliminary calculations with this. We also reduce the proofs of Key Propositions \ref{keyprop1} and \ref{keyprop2} to proving two analogous results (with more explicit dependence on the parameters $\mathcal{T},P,\beta$), which we call Key Propositions \ref{keyprop3} and \ref{keyprop4}, on the random multiplicative function side. The lower bound proposition, Key Proposition \ref{keyprop3}, is then proved in section \ref{seccorrlower}, and the upper bounds in Key Proposition \ref{keyprop4} are proved in section \ref{secfourthupper}.

In this paper, it isn't essential to work with random multiplicative functions. At the cost of a little technicality in places (one couldn't work with Euler products, but would always need approximations to them), and setting up appropriate initial analogues of Probability Results \ref{probres4} and \ref{probres5}, all of the proofs could be rewritten purely using character sums, Dirichlet polynomials and their averages. Some readers may be more comfortable with such a formulation. However, the author is strongly in favour of passing to random multiplicative functions (using \eqref{eqnorthcharrmf}) as quickly as possible, both to streamline the arguments and to emphasise their conceptual origins.

\subsection{Further remarks}
We close with some comparisons between our arguments here, and related results in the literature. Three obvious points of comparison are the proofs of analogous bounds for moments $\E|\sum_{n \leq x} f(n)|^{2q}$ of random multiplicative functions; the proofs of the upper bounds \eqref{eqbtscupper} and \eqref{eqbtsccontupper}; and other mollification-style or resonance-style arguments.

The proof~\cite{harperrmflow} of the bound $\E|\sum_{n \leq x} f(n)|^{2q} \asymp (\frac{x}{1 + (1-q)\sqrt{\log\log(10x)}} )^q$ (for $0 \leq q \leq 1$) is ostensibly quite unlike the proofs of our Theorems \ref{thmlb1} to \ref{thmlb3}. The upper and lower bound arguments are slightly different, but both begin by {\em conditioning} on the behaviour of $f(p)$ on many primes $p$ (all $(f(p))_{p \leq \sqrt{x}}$, in the lower bound case), before an application of the conditional H\"older inequality or the conditional Khintchine inequality. As discussed in some detail in the author's survey paper~\cite{harperbtsqc}, one cannot successfully perform anything like a conditioning on all $(\chi(p))_{p \leq \sqrt{x}}$, unless perhaps $x$ is very small compared with the conductor $r$ (e.g. of size $\log^{O(1)}r$). The starting point for the proof given here was the desire for an alternative, more ``low tech'' way to run this argument. In the situation of relevance, the proof of the conditional Khintchine inequality amounts to a comparison of conditional second and fourth moments. Key Propositions \ref{keyprop1} and \ref{keyprop2} are a substitute for this which can be executed for character sums (since they only depend on correlation calculations with character sums/Dirichlet polynomials), with the ``barrier'' inside our proxy object $I(\chi)$ taking the place of any global conditioning. Although the details of the calculations look rather different than in \cite{harperrmflow}, one ultimately arrives at square averages of Euler products that are close to those arising there. The final (and perhaps most interesting) part of the argument, exploiting the barrier to understand these Euler product averages, is actually then very similar as in the random case~\cite{harperrmflow} (see the discussion at the end of section \ref{subsecproofideas} immediately above, and especially section \ref{subsecexploitbarr}, below). Some of our contour integral manipulations here, especially in section \ref{subsecinitcond}, also have quite a lot in common with the preliminary conditional covariance calculations of Harper~\cite{harperas}, and conditional variance calculations of S. Hardy~\cite{shardy}.

In the author's paper~\cite{harpertypicalchar}, where the upper bounds \eqref{eqbtscupper} and \eqref{eqbtsccontupper} were proved, the same issue arose of requiring a substitute for the initial conditioning procedure used in the random case. However, in the upper bound argument there is some extra flexibility as compared with the lower bound argument, in that it suffices to condition on all $(f(p))_{p \leq P}$ for any parameter $P$, not too close to $x$ and satisfying $\log\log P \asymp \log\log x$. Notice one can choose such $P$ {\em much} smaller than $x$. In \cite{harpertypicalchar}, an analogue of this conditioning is developed for character sums, that works provided $e^{\log^{O(1)}P} < r/x$. Choosing $P = e^{\log^{c}L}$ (recall that $L = L_r := \min\{x, r/x\}$), we have $\log\log P \asymp \log\log L$, giving us the saving $\log\log L$ (in place of $\log\log x$) in \eqref{eqbtscupper}. And, as discussed previously, this is in fact the best saving we can hope for in general when working with character sums. We are permitted to work with such small $P$ in the upper bound argument because it relies on understanding ``conditioned'' second moments, and primes larger than $P$ (excluded from the ``conditioning'') simply make a mean square contribution corresponding to squareroot cancellation, but not worse than that. In contrast, lower bounds require working with something like conditioned fourth moments as well, and in the (more rapidly growing) fourth moment the primes larger than $P$ produce a blow-up of the shape $(\frac{\log x}{\log P})^{O(1)}$. Accommodating this would require keeping $P$ much closer to $x$, in general violating the condition $e^{\log^{O(1)}P} < r/x$ (the reader may compare with our discussion of the $P$-rough contribution and the respective sizes of $P$ and $x^{2\beta}$ in our argument here).

Conceptually, the key point is that in the conditioning analogue built in \cite{harpertypicalchar}, one must keep simultaneous control of prime number sums of lengths up to $P$ (or Euler products) at $\approx \log P$ different points, generating character sums of length $P^{\log^{O(1)}P}$. In our arguments here, the barrier $G_{\chi}(s)$ is only analysed at one or two points $s$ at a time, at $\approx \log\log P$ different scales (i.e. different truncation lengths $P^{e^{-j}}$), generating character sums of length $P^{(\log\log P)^{O(1)}}$. Thus we arrive (when attempting to invoke \eqref{eqnorthcharrmf}) at a much weaker condition like $P^{(\log\log P)^{O(1)}} < r/x$, permitting much larger $P$ for which our blow-up $(\frac{\log(x^{2\beta})}{\log P})^{O(1)}$ remains manageable.

Finally, consider some of the general landscape of mollification, resonance, and other such arguments. One could interpret our ``barrier'' $G_{\chi}(s)$ as something like a mollifier, being a Dirichlet polynomial/character sum that suppresses ``abnormal'' contributions from $\sum_{\substack{n \leq x^{\beta}, \\ P \; \text{smooth}}} \frac{\chi(n)}{n^s}$. The biggest difference from the classical set-up is perhaps that $G_{\chi}(s)$ is sensitive to the behaviour of all partial sums (or products) over $P^{e^{-j}}$-smooth numbers, not simply the full sum $\sum_{\substack{n \leq x^{\beta}, \\ P \; \text{smooth}}} \frac{\chi(n)}{n^s}$. Another important issue is that for our investigation of $\sum_{n \leq x} \chi(n)$, we crucially need to work with Dirichlet polynomials at a range of points $s$ inside the integral $I(\chi)$. In the more usual situation of trying to mollify $L(1/2, \chi)$, say, (or when investigating its moments, etc.), one would construct the mollifier as a Dirichlet polynomial at the single point $1/2$.

The argument of La Bret\`eche, Munsch and Tenenbaum~\cite{bretechemunschten}, delivering the best previously known lower bounds for $\frac{1}{r-1} \sum_{\chi \; \text{mod} \; r} |\sum_{n \leq x} \chi(n)|^{2q}$, $\frac{1}{r-1} \sum_{\chi \; \text{mod} \; r} |\sum_{n \leq x} h(n) \chi(n)|^{2q}$ and $\frac{1}{T} \int_{0}^{T} |\sum_{n \leq x} n^{it}|^{2q} dt$, also proceeds by introducing a proxy character sum and comparing various moments. In that argument, the proxy sum is constrained to have real non-negative coefficients. La Bret\`eche, Munsch and Tenenbaum~\cite{bretechemunschten} show that a choice of coefficients $c_n$ based on the ``anatomy'' of $n$ (i.e. the number of prime factors of $n$ of various sizes) delivers lower bounds that are off from the truth by powers of $\log L$, and also that there is {\em no choice} of real non-negative coefficients that can deliver substantially better bounds in their argument. Note that our proxy $I(\chi)$, if organised as a character sum, has complex coefficients. We also refer the reader to Szab\'o's paper~\cite{szabolower} proving lower bounds for the high moments $\frac{1}{r-1} \sum_{\substack{\chi \; \text{mod} \; r, \\ \chi \neq \chi_{0}}} |\sum_{n \leq x} \chi(n)|^{2q}$, where $q \geq 2$. This uses the same strategy, and a proxy object $R(\chi)$ that is quite closely related to our $I(\chi)$, namely a (discrete approximation to an) integral constructed to mimic appropriate properties of $|\sum_{n \leq x} \chi(n)|^{2(q-1)}$. However, because $q \geq 2$ some important features of our argument are not required in Szab\'o's case~\cite{szabolower}. In particular, no barrier $G_{\chi}(s)$ is needed, and the analogue of Key Proposition \ref{keyprop1} is much less delicate because $|\sum_{n \leq x} \chi(n)|$ can be replaced with $|\sum_{n \leq x} \chi(n)|^2$, so one can work throughout with obviously non-negative quantities rather than passing through $\Echar \sum_{n \leq x} \chi(n) \overline{I(\chi)}$.

\vspace{12pt}
Let us end by pointing out that between the author's previous paper~\cite{harpertypicalchar}, our results here, and Szab\'o's papers~\cite{szaboupper,szabolower}, we have a fairly complete understanding of the moments of (unweighted) character sums up to the second moment, and from the fourth moment onwards. This leaves a rather provocative gap, namely understanding the behaviour of $\frac{1}{r-1} \sum_{\chi \neq \chi_0 \; \text{mod} \; r} |\sum_{n \leq x} \chi(n)|^{2q}$ for $1 < q < 2$. On the random multiplicative function side, these moments are understood thanks to work of the author~\cite{harperrmfhigh}, and it seems reasonable to conjecture that (as in Szab\'o's~\cite{szaboupper,szabolower} high moment bounds) they should be $\asymp x^q \log^{(q-1)^2}(10L)$. We leave this as a challenge to the reader.

\section{Tools}\label{sectools}

\subsection{Smooth approximating functions}
In order to insert our required ``barrier'' conditions in a way that we can subsequently analyse, we shall use smooth functions that approximate (in fact, that closely upper bound) indicator functions. 

\begin{approxres1}\label{apres1}
For any $R \geq 0$ and $\delta > 0$, there exists a function $\gamma : \R \rightarrow [0,1+\delta]$ with the following properties:
\begin{enumerate}
\item $\gamma(x) \geq 1$ for all $|x| \leq R$;

\item $\gamma(x) \leq \delta$ for all $|x| > R + 1$;

\item for all $l \in \N$ and all $x \in \R$, we have the derivative estimate $|\frac{d^{l}}{dx^{l}} \gamma(x)| \leq \frac{(2R+1) (1+\delta)}{\pi (l + 1)} (\frac{2\pi}{\delta})^{l+1}$.

\end{enumerate}
\end{approxres1}

\begin{proof}[Proof of Approximation Result \ref{apres1}]
Results of roughly this shape occur throughout analysis, this particular statement is Approximation Result 1 of Harper~\cite{harperpartition}.
\end{proof}

\subsection{Number theory and harmonic analysis}
We record some number theoretic and analytic results that will be needed at various points in our proofs.

We shall require a suitable mean value estimate for Dirichlet polynomials, to handle the tails in some Perron integrals. There is much room for manoeuvre here, and e.g. the following result is convenient and more than good enough for us.

\begin{numth1}\label{numth1}
Uniformly for any complex numbers $(a_n)_{n=1}^{\infty}$ and any $T \geq 1$, and with $\Lambda(n)$ denoting the von Mangoldt function, we have
$$ \int_{-T}^{T} \Biggl|\sum_{T^{1.01} \leq n \leq x} \frac{a_n \Lambda(n)}{n^{1+it}} \Biggr|^{2} dt \ll \sum_{T^{1.01} \leq n \leq x} \frac{|a_n|^2 \Lambda(n)}{n} . $$
\end{numth1}

\begin{proof}[Proof of Number Theory Result \ref{numth1}]
Again, this is a relatively standard and classical type of bound. See e.g. Number Theory Result 1 of Harper~\cite{harperas} for a (brief) proof of this statement.
\end{proof}

Good pointwise bounds for (logarithmically weighted) Dirichlet polynomials over primes will be needed frequently, both in themselves and as ingredients in bounding other Dirichlet polynomials.

\begin{numth2}\label{numth2}
For any $100 \leq x \leq y$ (say), and any $t \neq 0$, we have
$$ \left|\sum_{x < p \leq y} \frac{\log p}{p^{1+it}} \right| \leq \frac{2}{|t|} + O((1+|t|) e^{-c\sqrt{\log x}}) , $$
$$ \left|\sum_{x < p \leq y} \frac{1}{p^{1+it}} \right| \leq \frac{2}{|t|\log x} + O((1+|t|) e^{-c\sqrt{\log x}}) . $$
\end{numth2}

\begin{proof}[Proof of Number Theory Result \ref{numth2}]
This follows by applying partial summation to the Prime Number Theorem with classical error term.
\end{proof}

An immediate application of Number Theory Result \ref{numth2} is the proof of the following important estimate.

\begin{numth3}\label{numth3}
Uniformly for all $100 \leq P \leq z$ and $1/\log P \leq |h| \leq e^{c\sqrt{\log P}}$, say (where $c > 0$ is a suitable small absolute constant, not necessarily the same as in Number Theory Result \ref{numth2}), we have
$$ \sum_{\substack{n \leq z, \\ n \; \text{is} \; P \; \text{rough}}} \frac{1}{n^{1+ih}} \ll 1 . $$
If $|h| \leq 1/\log P$, then instead
$$ \sum_{\substack{n \leq z, \\ n \; \text{is} \; P \; \text{rough}}} \frac{1}{n^{1+ih}} \ll \min\{\frac{1}{|h|\log P}, \frac{\log z}{\log P}\} . $$
\end{numth3}

If $\sum_{\substack{n \leq z, \\ n \; \text{is} \; P \; \text{rough}}} \frac{1}{n^{1+ih}}$ behaved like the analogous Euler product $\prod_{P < p \leq z} (1 - \frac{1}{p^{1+ih}})^{-1} \approx \exp\{\sum_{P < p \leq z}  \frac{1}{p^{1+ih}}\}$, then Number Theory Result \ref{numth3} would be a more or less immediate consequence of Number Theory Result \ref{numth2}. In particular, one should expect a bound $\ll 1$ when $|h| \geq 1/\log P$ (and $h$ isn't extremely large) because on this range the phases $p^{-ih}$ oscillate non-trivially, and so one sees substantial cancellation in $\sum_{P < p \leq z}  \frac{1}{p^{1+ih}}$. The rigorous proof proceeds along similar lines.

\begin{proof}[Proof of Number Theory Result \ref{numth3}]
We begin with the first part. We can write
$$ \sum_{\substack{n \leq z, \\ n \; \text{is} \; P \; \text{rough}}} \frac{1}{n^{1+ih}} = \prod_{P < p \leq z} (1 - \frac{1}{p^{1+ih}})^{-1} - \sum_{\substack{n > z, \\ p|n \Rightarrow P < p \leq z}} \frac{1}{n^{1+ih}} , $$
and Number Theory Result \ref{numth2} implies that the product $\prod_{P < p \leq z} (1 - \frac{1}{p^{1+ih}})^{-1}$ is $= \exp\{\sum_{P < p \leq z}  \frac{1}{p^{1+ih}} + O(1)\} = \exp\{O(\frac{1}{|h|\log P} + 1)\}$. If $|h| \geq 1/\log P$, this indeed contributes $\ll 1$. Using Abel summation, the subtracted sum is
\begin{eqnarray}
& \ll & \frac{1}{\log z} \max_{Z \geq z} \Biggl| \sum_{\substack{z < n \leq Z, \\ p|n \Rightarrow P < p \leq z}} \frac{\log n}{n^{1+ih}} \Biggr| \nonumber \\
& = & \frac{1}{\log z} \max_{Z \geq z} \Biggl| \sum_{\substack{z < n \leq Z, \\ p|n \Rightarrow P < p \leq z}} \frac{\sum_{d|n} \Lambda(d)}{n^{1+ih}} \Biggr| = \frac{1}{\log z} \max_{Z \geq z} \Biggl| \sum_{\substack{m : \\ p|m \Rightarrow P < p \leq z}} \frac{1}{m^{1+ih}} \sum_{\substack{z/m < d \leq Z/m, \\ p|d \Rightarrow P < p \leq z}} \frac{\Lambda(d)}{d^{1+ih}} \Biggr| \nonumber
\end{eqnarray}
Separating primes from prime powers and again using Number Theory Result \ref{numth2}, each of the sums over $d$ here is $= \sum_{\max\{P, z/m\} < p \leq \min\{z,Z/m\}} \frac{\log p}{p^{1+ih}} + O(1/P) \ll \frac{1}{|h|}$, so we get an overall contribution $\ll \frac{1}{|h|\log z}  \sum_{\substack{m : \\ p|m \Rightarrow P < p \leq z}} \frac{1}{m} \ll \frac{1}{|h|\log z} \frac{\log z}{\log P} = \frac{1}{|h|\log P} \ll 1$.

The second part is trivial when $|h| \leq 1/\log z$. When $1/\log z \leq |h| \leq 1/\log P$, it may be proved exactly similarly to the first part, but estimating the product there by $\prod_{P < p \leq e^{1/|h|}} (1 - \frac{1}{p^{1+ih}})^{-1} \cdot \prod_{e^{1/|h|} < p \leq z} (1 - \frac{1}{p^{1+ih}})^{-1} \ll \prod_{P < p \leq e^{1/|h|}} (1 - \frac{1}{p})^{-1} \ll \frac{1}{|h|\log P}$.
\end{proof}

At a couple of points in the proofs of Key Propositions \ref{keyprop3} and \ref{keyprop4} below (and thus ultimately of Key Propositions \ref{keyprop1} and \ref{keyprop2}), introducing appropriate smoothing will greatly simplify matters. We shall employ the classical Fej\'er kernel, giving us useful positivity, as well as quadratic decay of exponential integrals (as $|h|$ becomes large) rather than the linear decay that we would have using sharp cutoffs.

\begin{harman1}[Fej\'er kernel]\label{harman1}
For any $R > 0$ and any $h \in \R$, we have
$$ \int_{-R}^{R} e^{-ihw} \left(1 - \frac{|w|}{R} \right) dw = \frac{4\sin^{2}(hR/2)}{h^2 R} \ll \min\{R, \frac{1}{h^2 R} \} . $$
Furthermore, if $|h| \leq 1/R$ then the integral is also $\gg R$.
\end{harman1}

\begin{proof}[Proof of Harmonic Analysis Result \ref{harman1}]
The evaluation of the integral is a computation from classical Fourier analysis (and may easily be performed directly, e.g. using integration by parts). The claimed upper and lower bounds follow from the fact that $|\sin x| \asymp |x|$ when $|x| \leq 1$ (say), and $|\sin x| \leq 1$ for all $x \in \R$. See e.g. chapters 46 and 49 of K\"orner~\cite{korner} for some general discussion around the Fej\'er kernel.
\end{proof}

\subsection{Probabilistic tools}
To handle various secondary terms that arise in our arguments, we shall use the following rough bound for high moments of random multiplicative functions.

\begin{probres1}[Rough hypercontractive inequality]\label{probres1}
For any real $q \geq 1$, the following is true.

If $f(n)$ is a Steinhaus random multiplicative function, then for any sequence of complex numbers $(a_{n})_{n \leq N}$ we have
$$ \E\left|\sum_{n \leq N} a_n f(n) \right|^{2q} \leq \left(\sum_{n \leq N} |a_n|^2 d_{\lceil q \rceil}(n) \right)^{q} , $$
where $d_{k}(\cdot)$ denotes the $k$-fold divisor function (i.e. the number of $k$-tuples of natural numbers whose product is $\cdot$, or equivalently the Dirichlet series coefficient of $\zeta(s)^k$), and $\lceil q \rceil$ denotes the ceiling of $q$.
\end{probres1}

\begin{proof}[Proof of Probability Result \ref{probres1}]
This is a fairly simple and standard estimate. It may be proved by using H\"older's inequality to pass to the case of integer $q$, and then carefully expanding the $2q$-th power. See e.g. Probability Result 1 of Harper~\cite{harperrmfhigh} for a full proof.
\end{proof}

For ease of reference, we also record the following special case of Probability Result \ref{probres1}. This is sharp, apart from the exact dependence of the constant multiplier on $q$. With $q$ large, we shall apply this to control the (average) size of the error term in certain Taylor expansions, and with $q=2$ it will be an important preliminary step (after some conditioning) in the proof of Key Proposition \ref{keyprop4} below.

\begin{probres2}[Special case of Khintchine's inequality]\label{probres2}
In the setting of Probability Result \ref{probres1}, we have
$$ \E\Biggl|\sum_{\substack{p \leq N, \\ p \; \text{prime}}} a_p f(p) \Biggr|^{2q} \leq \lceil q \rceil^{q} \Biggl(\sum_{\substack{p \leq N, \\ p \; \text{prime}}} |a_p|^2 \Biggr)^{q} . $$
\end{probres2}

\begin{proof}[Proof of Probability Result \ref{probres2}]
For $p$ prime we have $d_{\lceil q \rceil}(p) = \lceil q \rceil$, so the result follows immediately from Probability Result \ref{probres1}. (Khintchine's inequality gives such a bound for all $q \geq 0$, and gives lower bounds as well. See e.g. Lemma 3.8.1 of Gut~\cite{gut}.) 
\end{proof}

Later in the proof of Key Proposition \ref{keyprop4}, we shall use the following estimate to control the correlations of certain terms.

\begin{probres3}\label{probres3}
Let $f(n)$ denote a Steinhaus random multiplicative function. Uniformly for all $100 \leq P \leq z$ and $1/\log P \leq |h| \leq e^{c\sqrt{\log P}}$, say (where $c > 0$ is a suitable small absolute constant), we have
$$ \E\Biggl|\sum_{\substack{n \leq z, \\ n \; \text{is} \; P \; \text{rough}}} \frac{f(n)}{n^{1/2}} \Biggr|^2 \Biggl|\sum_{\substack{n \leq z, \\ n \; \text{is} \; P \; \text{rough}}} \frac{f(n)}{n^{1/2+ih}} \Biggr|^2 \ll (\frac{\log z}{\log P})^2 . $$
\end{probres3}

Note that $ \E|\sum_{\substack{n \leq z, \\ n \; \text{is} \; P \; \text{rough}}} \frac{f(n)}{n^{1/2}} |^2 = \E|\sum_{\substack{n \leq z, \\ n \; \text{is} \; P \; \text{rough}}} \frac{f(n)}{n^{1/2+ih}} |^2 = \sum_{\substack{n \leq z, \\ n \; \text{is} \; P \; \text{rough}}} \frac{1}{n} \asymp \frac{\log z}{\log P}$, so Probability Result \ref{probres3} asserts that if $1/\log P \leq |h| \leq e^{c\sqrt{\log P}}$ then the shifted sums $\sum_{\substack{n \leq z, \\ n \; \text{is} \; P \; \text{rough}}} \frac{f(n)}{n^{1/2}}, \sum_{\substack{n \leq z, \\ n \; \text{is} \; P \; \text{rough}}} \frac{f(n)}{n^{1/2+ih}}$ are (roughly) uncorrelated.

\begin{proof}[Proof of Probability Result \ref{probres3}]
The content of the proof is actually purely number theoretic. Using the multiplicativity of $f$, we can rewrite the left hand side as
$$ \E \Biggl|\sum_{\substack{N \leq z^2, \\ N \; \text{is} \; P \; \text{rough}}} \frac{f(N) \sum_{\substack{m,n \leq z : \\ mn = N}} n^{-ih}}{N^{1/2}} \Biggr|^2 . $$
Then the orthogonality of the random variables $f(N)$ implies this is all equal to
$$ \sum_{\substack{N \leq z^2, \\ N \; \text{is} \; P \; \text{rough}}} \frac{|\sum_{\substack{m,n \leq z : \\ mn = N}} n^{-ih}|^2}{N} . $$

To proceed further, if we expand the square and let $\textbf{1}$ denote the indicator function, we can rewrite our sum as
$$ \sum_{\substack{n_1, n_2 \leq z, \\ P \; \text{rough}}} n_{1}^{-ih} n_{2}^{ih} \sum_{\substack{N \leq z^2, \\ N \; \text{is} \; P \; \text{rough}}} \frac{1}{N} \textbf{1}_{n_1 | N} \textbf{1}_{n_2 | N} \textbf{1}_{N/n_1 \leq z} \textbf{1}_{N/n_2 \leq z} . $$
Splitting up the sum over $n_1, n_2$ according to their highest common factor $d$, we find this is all
\begin{eqnarray}
& = & \sum_{\substack{d \leq z, \\ P \; \text{rough}}}  \sum_{\substack{m_1, m_2 \leq z/d, \\ (m_1,m_2) = 1, \\ P \; \text{rough}}} m_{1}^{-ih} m_{2}^{ih} \sum_{\substack{N \leq z^2, \\ N \; \text{is} \; P \; \text{rough}}} \frac{1}{N} \textbf{1}_{dm_1 | N} \textbf{1}_{dm_2 | N} \textbf{1}_{N/dm_1 \leq z} \textbf{1}_{N/dm_2 \leq z} \nonumber \\
& = & \sum_{\substack{d \leq z, \\ P \; \text{rough}}} \frac{1}{d} \sum_{\substack{m_1, m_2 \leq z/d, \\ (m_1,m_2) = 1, \\ P \; \text{rough}}} \frac{1}{m_{1}^{1+ih}} \frac{1}{m_{2}^{1-ih}} \sum_{\substack{M \leq z/\max\{m_1, m_2\}, \\ M \; \text{is} \; P \; \text{rough}}} \frac{1}{M} \nonumber \\
& = & \sum_{\substack{d \leq z, \\ P \; \text{rough}}} \frac{1}{d} \sum_{\substack{M \leq z, \\ M \; \text{is} \; P \; \text{rough}}} \frac{1}{M} \sum_{\substack{m_1, m_2 \leq \min\{z/d, z/M\}, \\ (m_1,m_2) = 1, \\ P \; \text{rough}}} \frac{1}{m_{1}^{1+ih}} \frac{1}{m_{2}^{1-ih}} . \nonumber
\end{eqnarray}

Finally, using the M\"obius function to detect the condition that the highest common factor $(m_1,m_2) = 1$, we can rewrite
\begin{eqnarray}
\sum_{\substack{m_1, m_2 \leq \min\{z/d, z/M\}, \\ (m_1,m_2) = 1, \\ P \; \text{rough}}} \frac{1}{m_{1}^{1+ih}} \frac{1}{m_{2}^{1-ih}} & = & \sum_{\substack{m_1, m_2 \leq \min\{z/d, z/M\}, \\ P \; \text{rough}}} \left(\sum_{e| (m_1,m_2)} \mu(e) \right) \frac{1}{m_{1}^{1+ih}} \frac{1}{m_{2}^{1-ih}} \nonumber \\
& = & \sum_{\substack{e \leq \min\{z/d, z/M\}, \\ P \; \text{rough}}} \frac{\mu(e)}{e^2} \sum_{\substack{m_1, m_2 \leq \min\{z/de, z/Me\}, \\ P \; \text{rough}}} \frac{1}{m_{1}^{1+ih}} \frac{1}{m_{2}^{1-ih}} . \nonumber
\end{eqnarray}
Number Theory Result \ref{numth3} implies this is all $\ll 1$. And since $\sum_{\substack{d \leq z, \\ P \; \text{rough}}} \frac{1}{d} \leq \prod_{P < p \leq z} (1-\frac{1}{p})^{-1} \ll \frac{\log z}{\log P}$, also for $\sum_{\substack{M \leq z, \\ M \; \text{is} \; P \; \text{rough}}} \frac{1}{M}$, our desired bound follows.
\end{proof}

At the heart of our arguments, we need estimates for expectations of random Euler products restricted by certain ``barrier'' events. When handling second moment type quantities (i.e. for Key Proposition \ref{keyprop3} below, and $\E |I(f)|^2$ in Key Proposition \ref{keyprop4}), the following ``one point'' estimate is what we shall need.

\begin{probres4}\label{probres4}
Let $f(n)$ denote a Steinhaus random multiplicative function. For any $P \geq 2$, we have
$$ \E \prod_{p \leq P} |1 - \frac{f(p)}{\sqrt{p}}|^{-2} \asymp \log P . $$
And there is a large number $B$ such that, uniformly for any large $a$, any large $P$, and any function $h(y)$ satisfying $|h(y)| \leq 100\log y$ (say), we have
\begin{eqnarray}
&& \E \prod_{p \leq P} |1 - \frac{f(p)}{\sqrt{p}}|^{-2} \textbf{1}_{-a - B(\log\log P - j) \leq \sum_{p \leq P^{e^{-j}}} \log|1 - \frac{f(p)}{\sqrt{p}}|^{-1} \leq a + (\log\log P - j) + h(\log\log P - j) \; \forall \; 0 \leq j \leq \log\log P - 1} \nonumber \\
& \asymp & \min\{1, \frac{a}{\sqrt{\log\log P}} \} \log P . \nonumber
\end{eqnarray}
\end{probres4}

The ultimate source of the saving factor $\sqrt{\log\log P}$ in this bound is the number of steps $j$ in the ``random walk'' $\sum_{p \leq P^{e^{-j}}} \log|1 - \frac{f(p)}{\sqrt{p}}|^{-1}$ inside the barrier condition. See the author's paper~\cite{harperrmflow}, and the survey~\cite{harperrmf3}, for much more discussion of this type of estimate.

\begin{proof}[Proof of Probability Result \ref{probres4}]
The first statement is easily proved by direct calculation, since we have
$$ \E \prod_{p \leq P} |1 - \frac{f(p)}{\sqrt{p}}|^{-2} = \E \Biggl|\sum_{\substack{n=1, \\ n \; \text{is} \; P \; \text{smooth}}}^{\infty} \frac{f(n)}{\sqrt{n}} \Biggr|^{2} =  \sum_{\substack{n=1, \\ n \; \text{is} \; P \; \text{smooth}}}^{\infty} \frac{1}{n} = \prod_{p \leq P} (1 - \frac{1}{p})^{-1} \asymp \log P . $$

For the second statement, note that it will suffice to prove this with the given event replaced by the condition that
\begin{eqnarray}
-a - B(\log\log P - j) & \leq & \sum_{\exp\{ e^{B} \frac{\log P}{e^{\lfloor \log\log P \rfloor}} \} < p \leq P^{e^{-j}}} \log|1 - \frac{f(p)}{\sqrt{p}}|^{-1} \nonumber \\
& \leq & a + (\log\log P - j) + h(\log\log P - j) \;\;\; \forall \; 0 \leq j \leq \log\log P - (B+1) , \nonumber
\end{eqnarray}
since the sum over primes smaller than $\exp\{ e^{B} \frac{\log P}{e^{\lfloor \log\log P \rfloor}} \}$ is uniformly bounded (and in particular has absolute value $\leq a/2$, say, provided $a$ is large enough), so may be included or discarded without essentially altering anything. With this adjustment, the claimed estimate follows from Proposition 5 of Harper~\cite{harperrmflow} with the choices $\sigma = 0$, $t_j \equiv 0$ and $n = \lfloor \log\log x \rfloor - (B+1)$, after interpreting the notation used there (e.g. replacing $x^{1/e}$ by $P$). Note also that the condition on the function $h$ in Proposition 5 of Harper~\cite{harperrmflow} is a little more restrictive (namely $|h(y)| \leq 10\log y$) than we impose here, but on inspecting the proofs it is easy to see that any fixed constant can be taken in place of 10 (with the required lower bounds on $B,P,a$ possibly then increasing).
\end{proof}

When handling fourth moment type quantities, as in the main part of Key Proposition \ref{keyprop4}, we shall need a ``two point'' estimate dealing with an Euler product and its shift.

\begin{probres5}\label{probres5}
In the setting of Probability Result \ref{probres4}, for any large $Q \leq P$ and any $t \in \R$ we have
$$ \E \prod_{Q < p \leq P} |1 - \frac{f(p)}{p^{1/2}}|^{-2} \cdot \prod_{Q < p \leq P} |1 - \frac{f(p)}{p^{1/2 + it}}|^{-2} \asymp \exp\{ \sum_{Q < p \leq P} \frac{2+2\cos(t\log p)}{p} \} . $$
In particular, if $1/\log Q \leq |t| \leq e^{c\sqrt{\log Q}}$ then this is $\asymp (\frac{\log P}{\log Q})^2$.

Furthermore, if $K \in \N$ and $P^{e^{-K}}$ are large; if $a$ is large, and $h(y)$ is a function satisfying $|h(y)| \leq 100\log y$; and if $e^{B}/\sqrt{\log(P^{e^{-K}})} \leq |t| \leq 1$ (say); then
\begin{eqnarray}
&& \E \prod_{P^{e^{-K}} < p \leq P} |1 - \frac{f(p)}{p^{1/2}}|^{-2} \cdot \prod_{P^{e^{-K}} < p \leq P} |1 - \frac{f(p)}{p^{1/2 + it}}|^{-2} \cdot \nonumber \\
&& \cdot \textbf{1}_{-a - B(K - j) \leq \sum_{P^{e^{-K}} < p \leq P^{e^{-j}}} \log|1 - \frac{f(p)}{p^{1/2}}|^{-1}, \; \sum_{P^{e^{-K}} < p \leq P^{e^{-j}}} \log|1 - \frac{f(p)}{p^{1/2 + it}}|^{-1} \leq a + (K - j) + h(K - j) \; \forall \; 0 \leq j \leq K - 1} \nonumber \\
& \asymp & \min\left\{1, \frac{a}{\sqrt{K}} \right\}^2 e^{2K} . \nonumber
\end{eqnarray}
\end{probres5}

The crucial point here is that, under the stated conditions on $t$, we find $\E \prod_{Q < p \leq P} |1 - \frac{f(p)}{p^{1/2}}|^{-2} \cdot \prod_{Q < p \leq P} |1 - \frac{f(p)}{p^{1/2 + it}}|^{-2}$ is $\asymp (\frac{\log P}{\log Q})^2 \asymp (\E \prod_{Q < p \leq P} |1 - \frac{f(p)}{p^{1/2}}|^{-2}) \cdot (\E \prod_{Q < p \leq P} |1 - \frac{f(p)}{p^{1/2 + it}}|^{-2})$, and similarly the expectation involving the simultaneous barrier conditions on $\sum_{P^{e^{-K}} < p \leq P^{e^{-j}}} \log|1 - \frac{f(p)}{p^{1/2}}|^{-1}$ and $\sum_{P^{e^{-K}} < p \leq P^{e^{-j}}} \log|1 - \frac{f(p)}{p^{1/2 + it}}|^{-1}$ is found to be $\asymp (\E \prod_{P^{e^{-K}} < p \leq P} |1 - \frac{f(p)}{p^{1/2}}|^{-2} \textbf{1}_{-a - B(K - j) \leq \sum_{P^{e^{-K}} < p \leq P^{e^{-j}}} \log|1 - \frac{f(p)}{p^{1/2}}|^{-1} \leq a + (K - j) + h(K - j) \; \forall \; 0 \leq j \leq K - 1}) \cdot (\E \prod_{P^{e^{-K}} < p \leq P} |1 - \frac{f(p)}{p^{1/2 + it}}|^{-2} \textbf{1}_{-a - B(K - j) \leq \sum_{P^{e^{-K}} < p \leq P^{e^{-j}}} \log|1 - \frac{f(p)}{p^{1/2 + it}}|^{-1} \leq a + (K - j) + h(K - j) \; \forall \; 0 \leq j \leq K - 1} )$. In other words, provided the shift $t$ is large enough the Euler products $\prod_{Q < p \leq P} |1 - \frac{f(p)}{p^{1/2}}|^{-2}, \prod_{Q < p \leq P} |1 - \frac{f(p)}{p^{1/2 + it}}|^{-2}$ become (roughly) uncorrelated, as do the events involving $\sum_{P^{e^{-K}} < p \leq P^{e^{-j}}} \log|1 - \frac{f(p)}{p^{1/2}}|^{-1}$ and $\sum_{P^{e^{-K}} < p \leq P^{e^{-j}}} \log|1 - \frac{f(p)}{p^{1/2 + it}}|^{-1}$.

\begin{proof}[Proof of Probability Result \ref{probres5}]
The first statement (again a reasonably straightforward direct calculation) follows from Lemma 6 of Harper~\cite{harperrmflow}, with the choices $u=v=\sigma = 0$. The estimate $\asymp (\frac{\log P}{\log Q})^2$ follows because the classical Mertens estimate, together with Number Theory Result \ref{numth2}, imply that $\sum_{Q < p \leq P} \frac{2+2\cos(t\log p)}{p} = 2(\sum_{Q < p \leq P} \frac{1}{p} + \Re \sum_{Q < p \leq P} \frac{1}{p^{1+it}}) = 2(\log\log P - \log\log Q) + O(1)$ under the given conditions on $t$.

The final part may be deduced from Proposition 7 of Harper~\cite{harperrmflow}, after interpreting the notation used there (e.g. replacing $x^{1/e}$ by $P$) and taking $\sigma = 0$, $D = \lfloor \log\log P \rfloor + 1 - K$ and $n = K$. Note that our assumption that $|t| \geq e^{B}/\sqrt{\log(P^{e^{-K}})}$ ensures that the important condition $D \geq 2\log(1/|t|) + B + 1$ in Proposition 7 is satisfied. As above, Proposition 7 is stated with a more restrictive condition $|h(y)| \leq 10\log y$, but in fact any fixed constant can be taken in place of 10.
\end{proof}

\section{Building the barrier, and passing to random multiplicative functions}\label{secpasstormf}
Before embarking on our main proofs, we must specify our choice of the ``barrier'' function $G_{\chi}(s) = G_{\chi, q,P}(s)$ in the definition of $I(\chi)$. We shall do this in two stages.

For ease of reference, let $B$ be the constant from Probability Results \ref{probres4} and \ref{probres5}, and $A$ be a suitably large absolute constant, and for each $t \in \R$ let $L_{\chi}(t) = L_{\chi, q,P}(t)$ denote the ``event'' that
\begin{eqnarray}
&& -\min\{\sqrt{\log\log P}, \frac{1}{1-q}\} - A + 2 -B\log(\frac{\log P}{e^j}) \nonumber \\
& \leq & \sum_{p \leq P^{e^{-j}}} \log|1 - \frac{\chi(p)}{p^{1/2 + it}}|^{-1} \nonumber \\
& \leq & \min\{\sqrt{\log\log P}, \frac{1}{1-q}\} + A - 2 + \log(\frac{\log P}{e^j}) - 50\log\log(\frac{\log P}{e^j}) \; \forall \; 0 \leq j \leq \log\log P - 1 . \nonumber
\end{eqnarray}
Let $U_{\chi}(t)$ denote the analogous event, with the numbers $2, -2$ replaced by $-3, 3$ respectively in the lower and upper bounds imposed.

Firstly, we define
$$ W_{\chi}(s) = W_{\chi, q,P}(s) := \left( \prod_{0 \leq j \leq \log\log P - 1} \gamma_{j}\Biggl( \Re \sum_{p \leq P^{e^{-j}}} \sum_{k=1,2} \frac{\chi(p^k)}{kp^{ks}} \Biggr) \right)^{\lfloor \log\log P \rfloor} , $$
where $\gamma_j$ is the function from Approximation Result \ref{apres1} with the choices $\delta = \frac{1}{(\log\log P)^2}$ and $R = R_j := \frac{(B+1)}{2}(\log\log P - j) - 25\log\log(\frac{\log P}{e^j}) + \min\{\sqrt{\log\log P}, \frac{1}{1-q}\} + A$, and with its argument shifted by $\frac{(B-1)}{2}(\log\log P - j) + 25\log\log(\frac{\log P}{e^j})$. Note that these choices certainly ensure that
\begin{equation}\label{eqWcrude}
0 \leq W_{\chi}(1/2 + it) \leq \left(1 + \frac{1}{(\log\log P)^2} \right)^{(\log\log P)^2} \ll 1 ,
\end{equation}
uniformly in $\chi$ and $t$. More precisely, Taylor expansion implies $\Re \sum_{p \leq P^{e^{-j}}} \sum_{k=1,2} \frac{\chi(p^k)}{kp^{k(1/2 + it)}} = \sum_{p \leq P^{e^{-j}}} \log|1 - \frac{\chi(p)}{p^{1/2 + it}}|^{-1} - \Re \sum_{p \leq P^{e^{-j}}} \sum_{k \geq 3} \frac{\chi(p^k)}{kp^{k(1/2 + it)}}$, and the absolute value of the subtracted sums is $\leq \sum_{p} \frac{1}{3} \sum_{k \geq 3} \frac{1}{p^{k/2}} \leq \sum_{p} \sum_{l \geq 1} \frac{1}{p^{3l/2}} \leq \sum_{n \geq 2} \frac{1}{n^{3/2}} \leq 2$, thus (using properties (i) and (ii) from Approximation Result \ref{apres1}) we get
\begin{equation}\label{eqWsharper}
\textbf{1}_{L_{\chi}(t)} \leq W_{\chi}(1/2 + it) \ll \textbf{1}_{U_{\chi}(t)} + \delta^{\lfloor \log\log P \rfloor} \ll \textbf{1}_{U_{\chi}(t)} + \frac{1}{\log^{1000}P} ,
\end{equation}
say. We take the large power $\lfloor \log\log P \rfloor$ in the definition of $W$ precisely to produce the very small term $1/\log^{1000}P$ in the upper bound in \eqref{eqWsharper}, without needing to choose $\delta$ itself prohibitively small. The fact that this error term is suitably tiny will be important when handling fourth moment type quantities in section \ref{subsecexploitbarr}, below.

We shall choose $G_{\chi}(s)$ to be a short character sum that is (usually) close to $W_{\chi}(s)$. The fact that the smooth functions $\gamma_j$ have well behaved derivatives (depending on $\delta$), as in part (iii) of Approximation Result \ref{apres1}, will enable us to do this. More precisely, let
$$ G_{\chi}(s) := \left( \prod_{0 \leq j \leq \log\log P - 1} \tilde{\gamma_{j}}\Biggl( \Re \sum_{p \leq P^{e^{-j}}} \sum_{k=1,2} \frac{\chi(p^k)}{kp^{ks}} \Biggr) \right)^{\lfloor \log\log P \rfloor} , $$
where $\tilde{\gamma_{j}}$ is the degree $2S-1$ Taylor polynomial of $\gamma_{j}$ about zero. For definiteness, let us take $S := \lfloor (\log\log P)^{10} \rfloor$. (The appropriateness of this choice will become apparent later, see e.g. section \ref{subsecgtow}.) Here we proceed somewhat similarly as in Proposition 1 of Harper~\cite{harpertypicalchar}. Then $G_{\chi}(s)$ may be expanded as a sum of character values, and their conjugates (thanks to the real part), of total length $\leq (P^2)^{2S(\log\log P)^2} \leq P^{4(\log\log P)^{12}}$. If we temporarily set
$$ \Delta = \Delta(\chi,s) := \max_{0 \leq j \leq \log\log P - 1} \Biggl|\gamma_{j}\Biggl( \Re \sum_{\substack{p \leq P^{e^{-j}}, \\ k=1,2}} \frac{\chi(p^k)}{kp^{ks}} \Biggr) - \tilde{\gamma_{j}}\Biggl( \Re \sum_{\substack{p \leq P^{e^{-j}}, \\ k=1,2}} \frac{\chi(p^k)}{kp^{ks}} \Biggr) \Biggr| , $$
then (arguing a little crudely, using the fact that $|\gamma_j| \leq 1+\delta$) we see
\begin{eqnarray}
|G_{\chi}(s) - W_{\chi}(s)| & \leq & (1 + \delta + \Delta )^{\lfloor \log\log P \rfloor^2} - (1+\delta)^{\lfloor \log\log P \rfloor^2} \nonumber \\
& \leq & \lfloor \log\log P \rfloor^2 \Delta \cdot (1 + \delta + \Delta )^{\lfloor \log\log P \rfloor^2 - 1} \nonumber \\
& \leq & (\log\log P)^2 \Delta \cdot 2^{(\log\log P)^2} \left( (1+\delta)^{\lfloor \log\log P \rfloor^2 - 1} + \Delta^{\lfloor \log\log P \rfloor^2 - 1} \right) . \nonumber
\end{eqnarray}
Combining Taylor's theorem with part (iii) of Approximation Result \ref{apres1}, we find $\Delta \ll \frac{1}{(2S)!} \frac{\log\log P}{S \delta} (\frac{2\pi}{\delta})^{2S} \max_{j} |\sum_{\substack{p \leq P^{e^{-j}}, \\ k=1,2}} \frac{\chi(p^k)}{kp^{ks}}|^{2S}$. Recalling that $\frac{1}{(2S)!} \leq (\frac{e}{2S})^{2S}$, and noting our choice of $S$, it follows that for a suitable large constant $C$ we will have 
\begin{equation}\label{eqGWcompare}
|G_{\chi}(s) - W_{\chi}(s)| \leq (\frac{C}{\delta S})^{2S} \max_{j} \Biggl|\sum_{\substack{p \leq P^{e^{-j}}, \\ k=1,2}} \frac{\chi(p^k)}{kp^{ks}} \Biggr|^{2S} + (\frac{C}{\delta S})^{2S \lfloor \log\log P \rfloor^2} \max_{j} \Biggl|\sum_{\substack{p \leq P^{e^{-j}}, \\ k=1,2}} \frac{\chi(p^k)}{kp^{ks}} \Biggr|^{2S \lfloor \log\log P \rfloor^2} .
\end{equation}

\vspace{12pt}
We are now ready to set things up for the proofs of Key Propositions \ref{keyprop1} and \ref{keyprop2}, by passing from character averages to the random multiplicative function side.

We clearly have $\Echar |\sum_{n \leq x} \chi(n)| |I(\chi)| \geq |\Echar \sum_{n \leq x} \chi(n) \overline{I(\chi)}|$, and in fact we shall prove the lower bound in Key Proposition \ref{keyprop1} for $ |\Echar \sum_{n \leq x} \chi(n) \overline{I(\chi)}|$.

If $f(n)$ is a Steinhaus random multiplicative function, we have
$$ \Echar \chi(n) \overline{\chi(m)} = \textbf{1}_{n=m} = \E f(n) \overline{f(m)} \;\;\;\;\;\;\;\; \forall \; 1 \leq n, m < r . $$
Consequently, and noting the lengths of the various character sums inside $I(\chi)$, we see that provided $\max\{x P^{4(\log\log P)^{12}}, x^{1-\beta} x^{\beta} x^{2\beta} P^{4(\log\log P)^{12}} \} = x^{1+2\beta} P^{4(\log\log P)^{12}} < r$ we will have $\Echar \sum_{n \leq x} \chi(n) \overline{I(\chi)} = \E \sum_{n \leq x} f(n) \overline{I(f)}$, and provided $x^{1+2\beta} P^{8(\log\log P)^{12}} < r$ also $\Echar |I(\chi)|^2 = \E |I(f)|^2$. Here $I(f)$ is defined in the obvious way, namely identically to $I(\chi)$ but with $\chi$ replaced everywhere by $f$. Under the stronger condition that $(x^{1-\beta} x^{\beta} x^{2\beta})^2 P^{16(\log\log P)^{12}} = x^{2+4\beta} P^{16(\log\log P)^{12}} < r$, we will also have $\Echar |I(\chi)|^4 = \E |I(f)|^4$. This leads us to:

\begin{keyprop3}\label{keyprop3}
There exists a small absolute constant $c > 0$ such that, uniformly for all large $x$, all $0 \leq q \leq 1$, and all small $\beta > 0$, the following is true. Provided that $e^{\log^{0.01}x} \leq P \leq x^{\beta/2}$ and $\log^{3}x \leq \mathcal{T} \leq e^{c\sqrt{\log P}}$ (say), we have
$$ |\E \sum_{n \leq x} f(n) \overline{I(f)}| \geq \frac{c\beta x}{1 + (1-q)\sqrt{\log\log P}} + O\left(\beta x \frac{\log^{3}x}{e^{\beta(\log x)/\log P}} + x \frac{\log P}{\log x} \log\log x \right) . $$
\end{keyprop3}

\begin{keyprop4}\label{keyprop4}
Uniformly for all large $x$, all $0 \leq q \leq 1$, and all small $\beta > 0$ such that $e^{\log^{0.01}x} \leq P \leq x^{\beta}$ and $0 \leq \mathcal{T} \leq e^{c\sqrt{\log P}}$ (say), we have
\begin{eqnarray}
\E |I(f)|^4 & \ll & e^{2\min\{\sqrt{\log\log P}, \frac{1}{1-q}\}} \left( \frac{\beta x}{1 + (1-q)\sqrt{\log\log P}} \right)^2 + \nonumber \\
&& + e^{2\min\{\sqrt{\log\log P}, \frac{1}{1-q}\}} \frac{(\beta x)^2}{(\log\log P)^{99}} \left( \frac{\beta\log x}{\log P} \right)^2 + \frac{\beta^4 x^2 \log^{2}x}{e^{2\beta(\log x)/\log P}} , \nonumber
\end{eqnarray}
and also
$$ \E |I(f)|^2 \ll \frac{\beta x}{1 + (1-q)\sqrt{\log\log P}} + \frac{\beta x}{e^{\beta(\log x)/\log P}} . $$
\end{keyprop4}

For clarity, let us reiterate that if we wish to employ Key Propositions \ref{keyprop3} and \ref{keyprop4} together (as we must to possibly deduce our main theorems), then we must make the same choices of $\beta, \mathcal{T}, P$ in both. There is enormous flexibility in the choice of $\mathcal{T}$, whereas $P$ must be chosen somewhat smaller than $x^{\beta}$, but not too much so, in order to keep the various error terms in Key Propositions \ref{keyprop3} and \ref{keyprop4} under control. This is not a technicality, but reflects real features of the problem. We need $P$ sufficiently smaller than $x^{\beta}$ that the sum $\sum_{\substack{n \leq x^{\beta}, \\ n \; \text{is} \; P \; \text{smooth}}} \frac{\chi(n)}{n^s}$ inside $F_{\chi}^{*}(s)$ (or the random analogue) behaves like an Euler product over primes $\leq P$. And we need $P$ sufficiently close to $x^{\beta}$ that only having control of this product (through the ``barrier'' $G_{\chi}(s)$), and not of the contribution from primes $> P$, is sufficient to control all of our integrals. It will suffice, for our purposes, to take
\begin{equation}\label{eqselectTP}
\mathcal{T} = \log^{3}x , \;\;\;\;\;\;\;\; \text{and} \;\;\;\;\;\;\;\; P = x^{\beta/(\log\log x)^{25}} .
\end{equation}

\begin{proof}[Proof of Key Proposition \ref{keyprop1}, assuming Key Proposition \ref{keyprop3}]
The hypotheses of Key Proposition \ref{keyprop1}, and our choice of $P$, guarantee (as in the above discussion) that
$$ \Echar |\sum_{n \leq x} \chi(n)| |I(\chi)| \geq |\Echar \sum_{n \leq x} \chi(n) \overline{I(\chi)}| = |\E \sum_{n \leq x} f(n) \overline{I(f)}| . $$
And with the choices of $\mathcal{T},P$ in \eqref{eqselectTP}, and our assumption that $\beta \geq 1/\log^{0.1}x$, Key Proposition \ref{keyprop3} is applicable and the ``big Oh'' term there is $\ll \frac{\beta x}{(\log\log x)^{24}}$, which is negligible compared with the main term $\frac{c\beta x}{1 + (1-q)\sqrt{\log\log P}} \gg \frac{\beta x}{1 + (1-q)\sqrt{\log\log x}}$.
\end{proof}

\begin{proof}[Proof of Key Proposition \ref{keyprop2}, assuming Key Proposition \ref{keyprop4}]
Similarly as for Key Proposition \ref{keyprop1}, our hypotheses and choice of $P$ guarantee that
$$ \Echar |I(\chi)|^2 = \E |I(f)|^2 , \;\;\;\;\; \text{and} \;\;\;\;\; \Echar |I(\chi)|^4 = \E |I(f)|^4 . $$
Key Proposition \ref{keyprop4} is applicable to estimate these, and the first terms in the bounds for $\E |I(f)|^4 , \E |I(f)|^2$ there are of the shape claimed in Key Proposition \ref{keyprop2}, so we need only check that the other terms are at most as big. With our choices we get $\frac{1}{(\log\log P)^{99}} ( \frac{\beta\log x}{\log P} )^2 \ll \frac{1}{(\log\log x)^{49}}$, so the second term in the bound for $\E |I(f)|^4$ is (more than) good enough. And the terms with denominators involving $e^{\beta(\log x)/\log P} = e^{(\log\log x)^{25}}$ make tiny contributions.
\end{proof}

Over the course of the next two sections, we shall work to prove Key Propositions \ref{keyprop3} and \ref{keyprop4}.

\section{Proof of Key Proposition \ref{keyprop3}}\label{seccorrlower}

\subsection{Improving the weight function}\label{subsecgtow}
We begin with a simple (although slightly drawn out) argument to show that $\E \sum_{n \leq x} f(n) \overline{I(f)}$ is close to $\E \sum_{n \leq x} f(n) \overline{J(f)}$, where we let
$$ J(f) := \frac{\sqrt{x}}{2\pi} \int_{-\mathcal{T}}^{\mathcal{T}} \Biggl(\sum_{x^{1-3\beta} < p \leq x^{1-\beta}} \frac{f(p)}{p^{1/2+iv}} \Biggr) F_{f}^{*}(1/2 + iv) W_{f}(1/2 + iv) \frac{x^{iv}}{1/2 + iv} dv $$
with the weight function $W$ as constructed above. This will allow us to work with $J(f)$, and thus exploit the nice properties \eqref{eqWcrude} and \eqref{eqWsharper} of $W$ (as opposed to its less agreeable approximation $G$), in the bulk of the proof.

Simply applying the triangle inequality, we have
\begin{eqnarray}
&& |\E \sum_{n \leq x} f(n) \overline{I(f)} - \E \sum_{n \leq x} f(n) \overline{J(f)}| \nonumber \\
& \leq & \sqrt{x} \int_{-\mathcal{T}}^{\mathcal{T}} \E\Biggl| \sum_{n \leq x} f(n) \Biggr| \Biggl| \sum_{\substack{x^{1-3\beta} \\ < p \leq x^{1-\beta}}} \frac{f(p)}{p^{1/2+iv}} \Biggr| |F_{f}^{*}(\frac{1}{2} + iv)| |G_{f}(\frac{1}{2} + iv) - W_{f}(\frac{1}{2} + iv)| \frac{dv}{|1/2 - iv|} . \nonumber
\end{eqnarray}
The Cauchy--Schwarz inequality, and orthogonality of the values $f(n)$, imply that the expectation here is
\begin{eqnarray}
& \leq & \sqrt{\E\Biggl| \sum_{n \leq x} f(n) \Biggr|^2} \sqrt{\E\Biggl| \sum_{x^{1-3\beta} < p \leq x^{1-\beta}} \frac{f(p)}{p^{1/2+iv}} \Biggr|^2 |F_{f}^{*}(\frac{1}{2} + iv)|^2 |G_{f}(\frac{1}{2} + iv) - W_{f}(\frac{1}{2} + iv)|^2} \nonumber \\
& \leq & \sqrt{x} \sqrt{\E\Biggl| \sum_{x^{1-3\beta} < p \leq x^{1-\beta}} \frac{f(p)}{p^{1/2+iv}} \Biggr|^2 |F_{f}^{*}(\frac{1}{2} + iv)|^2 |G_{f}(\frac{1}{2} + iv) - W_{f}(\frac{1}{2} + iv)|^2} . \nonumber 
\end{eqnarray}
Furthermore, using the independence of $f$ on distinct primes; the definition of $F_{f}^{*}(1/2 + iv)$; and the fact that $G_{f}(\frac{1}{2} + iv), W_{f}(\frac{1}{2} + iv)$ only depend on the $(f(p))_{p \leq P}$; the expectation under the second squareroot factors as
\begin{eqnarray}
&& \E\Biggl| \sum_{\substack{x^{1-3\beta} \\ < p \leq x^{1-\beta}}} \frac{f(p)}{p^{1/2+iv}} \Biggr|^2 \cdot \E|\sum_{\substack{n \leq x^{2\beta}, \\ P \; \text{rough}}} \frac{f(n)}{n^{1/2+iv}} |^2 \cdot \E |\sum_{\substack{n \leq x^{\beta}, \\ P \; \text{smooth}}} \frac{f(n)}{n^{1/2+iv}}|^2 |G_{f}(\frac{1}{2} + iv) - W_{f}(\frac{1}{2} + iv)|^2 \nonumber \\
& = & \sum_{\substack{x^{1-3\beta} \\ < p \leq x^{1-\beta}}} \frac{1}{p} \cdot \sum_{\substack{n \leq x^{2\beta}, \\ P \; \text{rough}}} \frac{1}{n} \cdot \E |\sum_{\substack{n \leq x^{\beta}, \\ P \; \text{smooth}}} \frac{f(n)}{n^{1/2+iv}}|^2 |G_{f}(\frac{1}{2} + iv) - W_{f}(\frac{1}{2} + iv)|^2 \nonumber \\
& \ll & 1 \cdot \frac{\log x}{\log P} \cdot \E |\sum_{\substack{n \leq x^{\beta}, \\ n \; \text{is} \; P \; \text{smooth}}} \frac{f(n)}{n^{1/2+iv}}|^2 |G_{f}(\frac{1}{2} + iv) - W_{f}(\frac{1}{2} + iv)|^2 . \nonumber
\end{eqnarray}

Finally, collecting our bounds together and applying the Cauchy--Schwarz inequality once more we get
\begin{eqnarray}\label{eqfinalgotoW}
&& |\E \sum_{n \leq x} f(n) \overline{I(f)} - \E \sum_{n \leq x} f(n) \overline{J(f)}| \\
& \ll & x \sqrt{\frac{\log x}{\log P}} \int_{-\mathcal{T}}^{\mathcal{T}} \Biggl( \E |\sum_{\substack{n \leq x^{\beta}, \\ P \; \text{smooth}}} \frac{f(n)}{n^{1/2+iv}}|^4 \Biggr)^{1/4} \left(  \E |G_{f}(\frac{1}{2} + iv) - W_{f}(\frac{1}{2} + iv)|^4 \right)^{1/4} \frac{dv}{|1/2 - iv|} . \nonumber
\end{eqnarray}
Using Probability Result \ref{probres1}, we find $\E |\sum_{\substack{n \leq x^{\beta}, \\ n \; \text{is} \; P \; \text{smooth}}} \frac{f(n)}{n^{1/2+iv}}|^4 \leq ( \sum_{\substack{n \leq x^{\beta}, \\ n \; \text{is} \; P \; \text{smooth}}} \frac{d(n)}{n} )^2 \leq (\prod_{p \leq P}(1 - \frac{1}{p})^{-2} )^2 \ll \log^{4}P$. Using \eqref{eqGWcompare}, and then upper bounding the maxima over $j$ there by sums, we deduce that $\E |G_{f}(\frac{1}{2} + iv) - W_{f}(\frac{1}{2} + iv)|^4$ is
$$ \ll \sum_{0 \leq j \leq \log\log P - 1} \Biggl( (\frac{C}{\delta S})^{8S} \E\Biggl|\sum_{\substack{p \leq P^{e^{-j}}, \\ k=1,2}} \frac{f(p^k)}{kp^{k(1/2+iv)}} \Biggr|^{8S} + (\frac{C}{\delta S})^{8S \lfloor \log\log P \rfloor^2} \E\Biggl|\sum_{\substack{p \leq P^{e^{-j}}, \\ k=1,2}} \frac{f(p^k)}{kp^{k(1/2+iv)}} \Biggr|^{8S \lfloor \log\log P \rfloor^2} \Biggr) . $$
We can easily and efficiently bound the expectations here by separating the prime and prime square contributions, and then applying Probability Result \ref{probres2} (noting that $f(p^2) = f(p)^2$ are again independent random variables, uniformly distributed on the unit circle). For example, we get
$$ \E\Biggl|\sum_{\substack{p \leq P^{e^{-j}}, \\ k=1,2}} \frac{f(p^k)}{kp^{k(1/2+iv)}} \Biggr|^{8S} \leq 2^{8S} \Biggl( \E\Biggl|\sum_{p \leq P^{e^{-j}}} \frac{f(p)}{p^{1/2 + iv}} \Biggr|^{8S} + \E\Biggl|\sum_{p \leq P^{e^{-j}}} \frac{f(p)}{2p^{1 + 2iv}} \Biggr|^{8S} \Biggr) \ll (16S)^{4S} \left( \sum_{p \leq P^{e^{-j}}} \frac{1}{p} \right)^{4S} , $$
which is $\ll (16S)^{4S} \left( \log\log P + O(1) \right)^{4S}$ in view of the Mertens estimate for the sum over primes. Exactly similarly, it follows that $\E|\sum_{\substack{p \leq P^{e^{-j}}, \\ k=1,2}} \frac{f(p^k)}{kp^{k(1/2+iv)}} |^{8S \lfloor \log\log P \rfloor^2}$ is $\ll (16S\lfloor \log\log P \rfloor^2 )^{4S\lfloor \log\log P \rfloor^2 } \left( \log\log P + O(1) \right)^{4S \lfloor \log\log P \rfloor^2}$. Thus, for a suitable absolute constant $C$ (different than above), we find \eqref{eqfinalgotoW} is
$$ \ll x \sqrt{\frac{\log x}{\log P}} \log\mathcal{T} \log P \log\log P \left( (\frac{C \log\log P}{\delta^2 S})^{S} + (\frac{C (\log\log P)^3}{\delta^2 S})^{S \lfloor \log\log P \rfloor^2} \right) . $$
Recalling that $\log\mathcal{T} \leq \sqrt{\log P}$, and our choices $\delta = 1/(\log\log P)^2$ and $S := \lfloor (\log\log P)^{10} \rfloor$, we conclude that \eqref{eqfinalgotoW} is $\ll \frac{x \sqrt{\log x}}{(\log\log P)^{(\log\log P)^{10}}}$. Under our assumption that $P \geq e^{\log^{0.01}x}$, this bound is more than good enough for Key Proposition \ref{keyprop3}.

\subsection{Conditioning and other initial manipulations}\label{subsecinitcond}
Our goal now is to produce a lower bound for $|\E \sum_{n \leq x} f(n) \overline{J(f)}|$. To do this, notice first that splitting up according to the largest prime factor $p$ of $n$, we can write
\begin{eqnarray}
\sum_{n \leq x} f(n) & = & \sum_{x^{1-3\beta} < p \leq x^{1-\beta}} f(p) \sum_{m \leq x/p} f(m) + \nonumber \\
&& + \sum_{x^{1-\beta} < p \leq x} f(p) \sum_{m \leq x/p} f(m) + \sum_{\substack{n \leq x, \\ n \; \text{is} \; x^{1-3\beta} \; \text{smooth}}} f(n) . \nonumber
\end{eqnarray}
Then it is easy to see (e.g. by first conditioning on the values $(f(p))_{p \leq x^{1-3\beta}}$) that only the sum $\sum_{x^{1-3\beta} < p \leq x^{1-\beta}} f(p) \sum_{m \leq x/p} f(m)$ will contribute to $\E \sum_{n \leq x} f(n) \overline{J(f)}$, and in fact
$$ \E (\sum_{n \leq x} f(n)) \overline{J(f)} = \frac{\sqrt{x}}{2\pi} \int_{-\mathcal{T}}^{\mathcal{T}} \sum_{\substack{x^{1-3\beta} \\ < p \leq x^{1-\beta}}} \frac{1}{p^{1/2-iv}} \E (\sum_{m \leq x/p} f(m)) \overline{F_{f}^{*}(1/2 + iv)} \overline{W_{f}(1/2 + iv)} \frac{x^{-iv}}{1/2 - iv} dv . $$

Next, the truncated Perron formula (with a somewhat crude handling of the error terms, and assuming that $\beta$ is small enough) implies that each sum $\sum_{m \leq x/p} f(m)$ here is
$$ = \frac{1}{2\pi i} \int_{1/2-ix^{0.1}}^{1/2+ix^{0.1}} \sum_{m \leq x^{3\beta}} \frac{f(m)}{m^s} \frac{(x/p)^s}{s} ds + O(1) = \frac{\sqrt{x}}{2\pi} \int_{-x^{0.1}}^{x^{0.1}} \sum_{m \leq x^{3\beta}} \frac{f(m)}{m^{1/2+it}} \frac{x^{it}}{p^{1/2+it}} \frac{dt}{1/2 + it} + O(1) . $$
See e.g. Corollary 5.3 of Montgomery and Vaughan~\cite{mv}. The contribution from all of these $O(1)$ terms to $\E (\sum_{n \leq x} f(n)) \overline{J(f)}$ is
$$ \ll \sqrt{x} \int_{-\mathcal{T}}^{\mathcal{T}} \sum_{x^{1-3\beta} < p \leq x^{1-\beta}} \frac{1}{p^{1/2}} \E |F_{f}^{*}(1/2 + iv)| \frac{1}{|1/2 - iv|} dv  , $$
where we used the fact that $W_f$ is uniformly bounded, as in \eqref{eqWcrude}. On using the easy bounds $\sum_{x^{1-3\beta} < p \leq x^{1-\beta}} \frac{1}{p^{1/2}} \leq \sum_{p \leq x^{1-\beta}} \frac{1}{p^{1/2}} \ll \frac{\sqrt{x^{1-\beta}}}{\log x}$ and $\E |F_{f}^{*}(1/2 + iv)| \leq \sqrt{\E |F_{f}^{*}(1/2 + iv)|^2} = \sqrt{(\sum_{\substack{n \leq x^{\beta}, \\ n \; \text{is} \; P \; \text{smooth}}} \frac{1}{n} ) \cdot (\sum_{\substack{n \leq x^{2\beta}, \\ n \; \text{is} \; P \; \text{rough}}} \frac{1}{n} )} \ll \sqrt{\beta \log x}$, we see this contribution is
$$ \ll \frac{x^{1-\beta/2} \sqrt{\beta\log x}}{\log x} \int_{-\mathcal{T}}^{\mathcal{T}} \frac{1}{|1/2 - iv|} dv \ll x^{1-\beta/2} \sqrt{\beta \log x} . $$
Meanwhile, the ``main term'' in our expression for $\E (\sum_{n \leq x} f(n)) \overline{J(f)}$ becomes
$$ \frac{x}{(2\pi)^2} \int_{-x^{0.1}}^{x^{0.1}} \int_{-\mathcal{T}}^{\mathcal{T}} \Biggl(\sum_{\substack{x^{1-3\beta} \\ < p \leq x^{1-\beta}}} \frac{1}{p^{1+i(t-v)}} \Biggr) \E \sum_{m \leq x^{3\beta}} \frac{f(m)}{m^{1/2+it}} \overline{F_{f}^{*}(\frac{1}{2} + iv)} \overline{W_{f}(\frac{1}{2} + iv)} x^{i(t-v)} \frac{dv}{1/2 - iv} \frac{dt}{1/2 + it} . $$

By conditioning now on the values $(f(p))_{p \leq P}$ that determine $\sum_{\substack{m \leq x^{\beta}, \\ P \; \text{smooth}}} \frac{f(m)}{m^{1/2+iv}}$ (inside $F_{f}^{*}(1/2 + iv)$) and $W_{f}(1/2 + iv)$, we may calculate that
\begin{eqnarray}\label{expcorrcalc}
&& \E \sum_{m \leq x^{3\beta}} \frac{f(m)}{m^{1/2+it}} \overline{F_{f}^{*}(\frac{1}{2} + iv)} \overline{W_{f}(\frac{1}{2} + iv)} \nonumber \\
& = & \E \Biggl( \sum_{\substack{n \leq x^{3\beta}, \\ P \; \text{rough}}} \frac{f(n)}{n^{1/2+it}} \sum_{\substack{m \leq x^{3\beta}/n, \\ P \; \text{smooth}}} \frac{f(m)}{m^{1/2+it}} \Biggr) \overline{\Biggl( \sum_{\substack{n \leq x^{2\beta}, \\ P \; \text{rough}}} \frac{f(n)}{n^{1/2+iv}} \sum_{\substack{m \leq x^{\beta}, \\ P \; \text{smooth}}} \frac{f(m)}{m^{1/2+iv}} \Biggr)} \overline{W_{f}(\frac{1}{2} + iv)} \nonumber \\
& = & \sum_{\substack{n \leq x^{2\beta}, \\ n \; \text{is} \; P \; \text{rough}}} \frac{1}{n^{1+i(t-v)}} \E  (\sum_{\substack{m \leq x^{3\beta}/n, \\ P \; \text{smooth}}} \frac{f(m)}{m^{1/2+it}}) \cdot (\sum_{\substack{m \leq x^{\beta}, \\ P \; \text{smooth}}} \frac{\overline{f(m)}}{m^{1/2 - iv}}) \overline{W_{f}(\frac{1}{2} + iv)} .
\end{eqnarray}
It will be convenient for our later calculations to now replace $(\sum_{\substack{m \leq x^{3\beta}/n, \\ m \; \text{is} \; P \; \text{smooth}}} \frac{f(m)}{m^{1/2+it}}) \cdot (\sum_{\substack{m \leq x^{\beta}, \\ m \; \text{is} \; P \; \text{smooth}}} \frac{\overline{f(m)}}{m^{1/2 - iv}})$ by $F_{P}(1/2+it) \overline{F_{P}(1/2+iv)}$, where $F_{P}(s) := \prod_{p \leq P} (1 - \frac{f(p)}{p^s})^{-1} = \sum_{\substack{m=1, \\ m \; \text{is} \; P \; \text{smooth}}}^{\infty} \frac{f(m)}{m^{s}}$ is the usual random Euler product corresponding to $f$ on $P$-smooth numbers. (This will remove the dependence on $n$ from these terms, as well as introducing the useful product structure.) We will be able to do this at the cost of an error term that is small provided $P$ is small enough compared with $x^{\beta}$. Indeed, since $W_f$ is uniformly bounded, making this change creates an error term inside the expectation in \eqref{expcorrcalc} that is
$$ \ll \sqrt{\E|\sum_{\substack{m \leq x^{3\beta}/n, \\ P \; \text{smooth}}} \frac{f(m)}{m^{1/2 + it}}|^2} \sqrt{\E|\sum_{\substack{m > x^{\beta}, \\ P \; \text{smooth}}} \frac{\overline{f(m)}}{m^{1/2 - iv}}|^2} + \sqrt{\E|F_{P}(1/2+iv)|^2} \sqrt{\E|\sum_{\substack{m > x^{3\beta}/n, \\ P \; \text{smooth}}} \frac{f(m)}{m^{1/2 + it}}|^2} . $$
Noting that $x^{3\beta}/n \geq x^{\beta}$ here, and using Rankin's trick, we may bound this error by
$$ \ll \sqrt{\log P} \sqrt{\sum_{\substack{m > x^{\beta}, \\ m \; \text{is} \; P \; \text{smooth}}} \frac{1}{m}} \ll \sqrt{\log P} \sqrt{\frac{1}{x^{2\beta/\log P}} \sum_{\substack{m: \\ m \; \text{is} \; P \; \text{smooth}}} \frac{1}{m^{1-2/\log P}}} \ll \frac{\log P}{e^{\beta(\log x)/\log P}} . $$
This produces an overall contribution to $\E (\sum_{n \leq x} f(n)) \overline{J(f)}$ that is
$$ \ll x \int_{-x^{0.1}}^{x^{0.1}} \int_{-\mathcal{T}}^{\mathcal{T}} \Biggl(\sum_{x^{1-3\beta} < p \leq x^{1-\beta}} \frac{1}{p} \Biggr) \sum_{\substack{n \leq x^{2\beta}, \\ n \; \text{is} \; P \; \text{rough}}} \frac{1}{n} \frac{\log P}{e^{\beta(\log x)/\log P}} \frac{dv}{|1/2 - iv|} \frac{dt}{|1/2 + it|} \ll x \frac{\beta \log^{3}x}{e^{\beta(\log x)/\log P}} . $$

We have so far shown that $\E (\sum_{n \leq x} f(n)) \overline{J(f)}$ is
\begin{eqnarray}\label{zerointred}
& = & \frac{x}{(2\pi)^2} \int_{-x^{0.1}}^{x^{0.1}} \int_{-\mathcal{T}}^{\mathcal{T}} \kappa(t-v) \E F_{P}(\frac{1}{2}+it) \overline{F_{P}(\frac{1}{2} + iv)} \overline{W_{f}(\frac{1}{2} + iv)} \frac{dv}{1/2 - iv} \frac{dt}{1/2 + it} + \nonumber \\
&& + O\left( x \frac{\beta \log^{3}x}{e^{\beta(\log x)/\log P}} \right) ,
\end{eqnarray}
where $\kappa(h) := (\sum_{x^{1-3\beta} < p \leq x^{1-\beta}} \frac{1}{p^{1+ih}}) \cdot (\sum_{\substack{n \leq x^{2\beta}, \\ n \; \text{is} \; P \; \text{rough}}} \frac{1}{n^{1+ih}}) \cdot x^{ih}$. For our final basic reduction, we will show that the ranges of integration in \eqref{zerointred} may be significantly restricted. This may be done very similarly as in section 3.2 of Harper~\cite{harperas}. Indeed, if we simply use the crude bounds $|\kappa(h)| \ll \frac{\log x}{\log P} |\sum_{x^{1-3\beta} < p \leq x^{1-\beta}} \frac{1}{p^{1+ih}}|$ and $|\E F_{P}(\frac{1}{2}+it) \overline{F_{P}(\frac{1}{2} + iv)} \overline{W_{f}(\frac{1}{2} + iv)}| \ll \E|F_{P}(1/2)|^2 \ll \log P$, it follows that the overall contribution from $\mathcal{T} < |t| \leq x^{0.1}$ is
\begin{eqnarray}
& \ll & x\log x \int_{-\mathcal{T}}^{\mathcal{T}} \frac{1}{|1/2 - iv|} \int_{\mathcal{T} < |t| \leq x^{0.1}} \frac{1}{|1/2 + it|} \Biggl| \sum_{x^{1-3\beta} < p \leq x^{1-\beta}} \frac{1}{p^{1+i(t-v)}} \Biggr| dt dv \nonumber \\
& \leq & x\log x \int_{-\mathcal{T}}^{\mathcal{T}} \frac{1}{|1/2 - iv|} \sqrt{\int_{\mathcal{T} < |t| \leq x^{0.1}} \frac{dt}{|1/2 + it|^2}} \sqrt{\int_{\mathcal{T} < |t| \leq x^{0.1}} \Biggl| \sum_{x^{1-3\beta} < p \leq x^{1-\beta}} \frac{1}{p^{1+i(t-v)}} \Biggr|^2 dt} dv . \nonumber
\end{eqnarray}
Under the squareroots here, the first integral over $t$ is clearly $\ll 1/\mathcal{T}$, and Number Theory Result \ref{numth1} implies the second integral over $t$ (in fact even when extended over all $|t| \leq x^{0.1}$) is $\ll \sum_{x^{1-3\beta} < p \leq x^{1-\beta}} \frac{1}{p\log p} \ll 1/\log x$. Thus we get a bound $\ll x (\log\mathcal{T}) \sqrt{\frac{\log x}{\mathcal{T}}}$ for the contribution from $|t| > \mathcal{T}$.

Since we assume that $\mathcal{T} \leq e^{c\sqrt{\log P}} \leq e^{c\sqrt{\log x}}$, Number Theory Result \ref{numth2} implies that $|\sum_{x^{1-3\beta} < p \leq x^{1-\beta}} \frac{1}{p^{1+i(t-v)}} | \ll \frac{1}{|t-v|\log x}$ for all $|t|, |v| \leq \mathcal{T}$, and Number Theory Result \ref{numth3} implies that $|\sum_{\substack{n \leq x^{2\beta}, \\ n \; \text{is} \; P \; \text{rough}}} \frac{1}{n^{1+i(t-v)}}| \ll 1$ for all $1/\log P \leq |t-v| \leq 2\mathcal{T}$. Thus $|\kappa(t-v)| \ll \frac{1}{|t-v|\log x}$ for such $t,v$. Again combining this estimate with the easy bound $|\E F_{P}(\frac{1}{2}+it) \overline{F_{P}(\frac{1}{2} + iv)} \overline{W_{f}(\frac{1}{2} + iv)}| \ll \E|F_{P}(1/2)|^2 \ll \log P$, we see the contribution to \eqref{zerointred} from all $|t|, |v| \leq \mathcal{T}$ with $|t-v| \geq 1/\log P$ is
$$ \ll x \frac{\log P}{\log x} \int \int_{\substack{|t|, |v| \leq \mathcal{T}, \\ |t-v| \geq 1/\log P}} \frac{1}{|t-v|} \frac{dv}{|1/2 - iv|} \frac{dt}{|1/2 + it|} . $$
In particular, when $1/\log P \leq |t-v| \leq 1$ we get a contribution that is
$$ \ll x \frac{\log P}{\log x} \int_{-\mathcal{T}}^{\mathcal{T}} \frac{1}{1 + t^2} \int_{1/\log P \leq |t-v| \leq 1} \frac{1}{|t-v|} dv dt \ll x \frac{\log P}{\log x} \log\log P . $$
When $|t-v| \geq 1$ (and so $|t-v| \asymp (1 + |t-v|)$), using the symmetry of $t,v$ we get a contribution that is
\begin{eqnarray}
& \ll & x \frac{\log P}{\log x} \int_{-\mathcal{T}}^{\mathcal{T}} \frac{1}{1+|t|} \left(\int_{|v| < \frac{|t|}{2}} \frac{dv}{(1+|t|)(1+|v|)} + \int_{\frac{|t|}{2} \leq |v| \leq |t|} \frac{dv}{(1+|t-v|)(1+|t|)} \right) dt \nonumber \\
& \ll &  x \frac{\log P}{\log x} \int_{-\mathcal{T}}^{\mathcal{T}} \frac{1}{1+|t|} \frac{\log(2+|t|)}{1+|t|} dt \ll x \frac{\log P}{\log x} . \nonumber
\end{eqnarray}

\subsection{Reducing to a single integral}\label{subsecreducesingle}
Summarising our work so far, we have shown that $\E (\sum_{n \leq x} f(n)) \overline{J(f)}$ is
\begin{eqnarray}\label{firstintred}
& = & \frac{x}{(2\pi)^2} \int \int_{\substack{|v| \leq \mathcal{T}, \\ |t-v| \leq 1/\log P}} \kappa(t-v) \E F_{P}(\frac{1}{2}+it) \overline{F_{P}(\frac{1}{2} + iv)} \overline{W_f(\frac{1}{2} + iv)} \frac{dv}{1/2 - iv} \frac{dt}{1/2 + it} + \nonumber \\
&& + O\left( x \frac{\beta \log^{3}x}{e^{\beta(\log x)/\log P}} + x (\log\mathcal{T}) \sqrt{\frac{\log x}{\mathcal{T}}} + x \frac{\log P}{\log x} \log\log P \right) ,
\end{eqnarray}
where we recall $\kappa(h) = (\sum_{x^{1-3\beta} < p \leq x^{1-\beta}} \frac{1}{p^{1+ih}}) \cdot (\sum_{\substack{n \leq x^{2\beta}, \\ n \; \text{is} \; P \; \text{rough}}} \frac{1}{n^{1+ih}}) \cdot x^{ih}$.

Our next goal is to actually perform the integration over the short $t$ variable. Since the Euler product $F_{P}(\frac{1}{2}+it)$ only runs over primes that are $\leq P$ (so involves terms $p^{-it} = e^{-it\log p}$ with $p \leq P$), and in \eqref{firstintred} we have $|t-v| \leq 1/\log P$, we may expect that typically $F_{P}(\frac{1}{2}+it) \approx F_{P}(\frac{1}{2}+iv)$ there. The following lemma will help to make this precise.

\begin{lem1}\label{lem1}
In the above situation, for any $t,v \in \R$ we have
$$ \E F_{P}(\frac{1}{2}+it) \overline{F_{P}(\frac{1}{2} + iv)} \overline{W_f(\frac{1}{2} + iv)} = \E |F_{P}(\frac{1}{2} + iv)|^2 \overline{W_f(\frac{1}{2} + iv)} + O(|t-v|\log^{2}P) . $$
\end{lem1}

\begin{proof}[Proof of Lemma \ref{lem1}]
The difference between $\E F_{P}(\frac{1}{2}+it) \overline{F_{P}(\frac{1}{2} + iv)} \overline{W_f(\frac{1}{2} + iv)}$ and $\E |F_{P}(\frac{1}{2} + iv)|^2 \overline{W_f(\frac{1}{2} + iv)}$ is certainly $\ll \E |F_{P}(\frac{1}{2}+it) - F_{P}(\frac{1}{2}+iv)||F_{P}(\frac{1}{2} + iv)|$, so in view of the Cauchy--Schwarz inequality and the easy estimate $\E |F_{P}(\frac{1}{2} + iv)|^2 \ll \log P$, to prove Lemma \ref{lem1} it will suffice to show that
$$ \E |F_{P}(\frac{1}{2}+it) - F_{P}(\frac{1}{2} + iv)|^2 \ll |t-v|^2 \log^{3}P . $$

But using orthogonality of the $f(m)$, we may calculate explicitly that the left hand side is
$$ = \E |\sum_{\substack{m=1, \\ m \; \text{is} \; P \; \text{smooth}}}^{\infty} \frac{f(m)(m^{-it} - m^{-iv})}{\sqrt{m}}|^2 = \sum_{\substack{m=1, \\ P \; \text{smooth}}}^{\infty} \frac{|m^{-it} - m^{-iv}|^2}{m} \ll |t-v|^2 \sum_{\substack{m=1, \\ P \; \text{smooth}}}^{\infty} \frac{\log^{2}m}{m} . $$
Then the sum over $m$ may be estimated in various standard ways, for example it is
$$ = \log^{2}P \sum_{\substack{m=1, \\ m \; \text{is} \; P \; \text{smooth}}}^{\infty} \frac{((\log m)/\log P)^2}{m} \ll \log^{2}P \sum_{\substack{m=1, \\ m \; \text{is} \; P \; \text{smooth}}}^{\infty} \frac{1}{m^{1-1/\log P}} \ll \log^{3}P . $$
\end{proof}

We may also note that $\frac{1}{1/2 + it} = \frac{1+O(|t-v|)}{1/2 + iv}$. Substituting this all into \eqref{firstintred}, and using the bound $|\kappa(h)| \ll \min\{\frac{\log x}{\log P}, \frac{1}{|h|^2 \log x \log P}\}$ for $|h| \leq 1/\log P$ (which follows from Number Theory Results \ref{numth2} and \ref{numth3}, similarly as above), we find that the first line of \eqref{firstintred} is
\begin{eqnarray}
& = & \frac{x}{(2\pi)^2} \int_{-\mathcal{T}}^{\mathcal{T}} \int_{-1/\log P}^{1/\log P} \kappa(h) (\E |F_{P}(\frac{1}{2} + iv)|^2 \overline{W_f(\frac{1}{2} + iv)} + O(|h|\log^{2}P)) dh \frac{dv}{1/4 + v^2} \nonumber \\
& = & \frac{x}{(2\pi)^2} \left(\int_{-\mathcal{T}}^{\mathcal{T}} \frac{\E |F_{P}(\frac{1}{2} + iv)|^2 \overline{W_f(\frac{1}{2} + iv)}}{1/4 + v^2} dv \right) \left( \int_{-1/\log P}^{1/\log P} \kappa(h) dh \right)  + O(x \frac{\log P}{\log x} \log\log x) . \nonumber
\end{eqnarray}

To analyse $\int_{-1/\log P}^{1/\log P} \kappa(h) dh$, a fairly neat and clean approach is to introduce a Fej\'er kernel (which will give us some useful positivity). Thus this integral is
$$ = \int_{-1/\log P}^{1/\log P} \kappa(h) (1 - |h|\log P) dh + O\left(\log P \int_{-1/\log P}^{1/\log P} |h| |\kappa(h)| dh \right) , $$
and again using the bound $|\kappa(h)| \ll \min\{\frac{\log x}{\log P}, \frac{1}{|h|^2 \log x \log P}\}$ shows the ``big Oh'' term is $\ll \int_{-1/\log P}^{1/\log P} \min\{\frac{1}{|h| \log x}, |h|\log x \} dh \ll \frac{\log\log x}{\log x}$. This contributes a further $O(x \log P \frac{\log\log x}{\log x})$ to \eqref{firstintred}. Meanwhile, recalling the explicit definition of $\kappa(h)$ and writing $(x/np)^{ih} = e^{ih\log(x/np)}$, and using the properties of the Fej\'er kernel from Harmonic Analysis Result \ref{harman1}, we see
\begin{eqnarray}
\int_{-1/\log P}^{1/\log P} \kappa(h) (1 - |h|\log P) dh & = & \sum_{x^{1-3\beta} < p \leq x^{1-\beta}} \frac{1}{p} \sum_{\substack{n \leq x^{2\beta}, \\ n \; \text{is} \; P \; \text{rough}}} \frac{1}{n} \frac{4\log P \sin^2(\frac{\log(x/np)}{2\log P})}{\log^{2}(x/np)} \nonumber \\
& \gg & \frac{1}{\log P} \sum_{x^{1-3\beta} < p \leq x^{1-\beta}} \frac{1}{p} \sum_{\substack{n \leq x^{2\beta}, \\ n \; \text{is} \; P \; \text{rough}}} \frac{1}{n} \textbf{1}_{x/P \leq np \leq xP} \nonumber \\
& \geq & \frac{1}{\log P} \sum_{\substack{x^{\beta}P \leq n \leq x^{2\beta}, \\ n \; \text{is} \; P \; \text{rough}}} \frac{1}{n} \sum_{x/nP \leq p \leq xP/n} \frac{1}{p} . \nonumber
\end{eqnarray}
A standard Mertens estimate shows the sum over primes here is $= \log(\frac{\log(xP/n)}{\log(x/nP)}) + O(\frac{1}{\log x}) = \log(1 + \frac{2\log P}{\log(x/nP)}) + O(\frac{1}{\log x}) \gg \frac{\log P}{\log x}$. And classical estimates for the counting function of rough numbers (see e.g. Theorem 6.4 of Tenenbaum~\cite{tenenbaumintro}), together with our assumption that $P \leq x^{\beta/2}$, imply that $ \sum_{\substack{x^{\beta}P \leq n \leq x^{2\beta}, \\ n \; \text{is} \; P \; \text{rough}}} \frac{1}{n} \gg \frac{\log(x^{\beta}/P)}{\log P} \gg \frac{\beta\log x}{\log P}$, so we obtain that $\int_{-1/\log P}^{1/\log P} \kappa(h) (1 - |h|\log P) dh \gg \frac{\beta}{\log P}$.

\subsection{Conclusion}\label{subseclowercon}
We have now established that, for a certain absolute constant $c > 0$, we have
\begin{eqnarray}
|\E \sum_{n \leq x} f(n) \overline{J(f)}| & \geq & \frac{c\beta x}{\log P} \int_{-\mathcal{T}}^{\mathcal{T}} \frac{\E |F_{P}(\frac{1}{2} + iv)|^2 \overline{W_f(\frac{1}{2} + iv)}}{1/4 + v^2} dv + \nonumber \\
&& + O\left( x \frac{\beta\log^{3}x}{e^{\beta(\log x)/\log P}} + x (\log\mathcal{T}) \sqrt{\frac{\log x}{\mathcal{T}}} + x \frac{\log P}{\log x} \log\log x \right) . \nonumber
\end{eqnarray}
It is perhaps worth remarking that we haven't used too much about the nature of $W_f(1/2 + iv)$ thus far, only its boundedness and the fact that it only depends on the values of $f$ on primes $\leq P$.

Our assumption that $\mathcal{T} \geq \log^{3}x$ (which is also only used for simplification at this point) implies that the error term $x (\log\mathcal{T}) \sqrt{\frac{\log x}{\mathcal{T}}}$ is negligible compared with $x \frac{\log P}{\log x} \log\log x$, so may be dropped.

Thanks to translation invariance in law, the integral over $v$ on the first line is $= \E |F_{P}(\frac{1}{2})|^2 \overline{W_f(\frac{1}{2})} (\int_{-\mathcal{T}}^{\mathcal{T}} \frac{dv}{1/4 + v^2}) \gg \E |F_{P}(\frac{1}{2})|^2 \overline{W_f(\frac{1}{2})}$, which is $\geq \E |F_{P}(\frac{1}{2})|^2 \textbf{1}_{L_{f}(0)}$ thanks to \eqref{eqWsharper}. Finally, applying Probability Result \ref{probres4} with the choices $a = \min\{\sqrt{\log\log P}, \frac{1}{1-q}\} + A - 2$ and $h(y) = -50\log y$, we get $\E |F_{P}(\frac{1}{2})|^2 \textbf{1}_{L_{f}(0)} \gg \min\{1, \frac{a}{\sqrt{\log\log P}} \} \log P \gg \frac{\log P}{1 + (1-q)\sqrt{\log\log P}}$. Note that provided $A$ is fixed large enough, $a$ will be sufficiently large that Probability Result \ref{probres4} may legitimately be applied.

Key Proposition \ref{keyprop3} now follows, for an appropriately adjusted value of $c$.
\qed

\section{Proof of Key Proposition \ref{keyprop4}}\label{secfourthupper}
We concentrate on establishing the claimed upper bound for $\E|I(f)|^4$. The case of $\E|I(f)|^2$, which is much more straightforward, will be briefly discussed at the end.

\subsection{Khintchine's inequality, and passing to an Euler product average}\label{subsecfourthkh}
Rather than expanding $|I(f)|^4$ and working directly with a fourfold integral, we can simplify matters with a short preliminary conditioning step.

Thus let $\tilde{\E}$ denote expectation conditional on the values $(f(p))_{p \leq x^{1-3\beta}}$ (under which, in particular, all the terms $F_{f}^{*}(1/2 + iv),  G_{f}(1/2 + iv)$ become fixed). The Tower Property of conditional expectation implies that $\E|I(f)|^4 = \E\tilde{\E} |I(f)|^4$. Since we can rewrite
$$ I(f) = \frac{\sqrt{x}}{2\pi} \sum_{x^{1-3\beta} < p \leq x^{1-\beta}} \frac{f(p)}{p^{1/2}} \int_{-\mathcal{T}}^{\mathcal{T}} \frac{1}{p^{iv}} F_{f}^{*}(1/2 + iv) G_{f}(1/2 + iv) \frac{x^{iv}}{1/2 + iv} dv , $$
an application of Probability Result \ref{probres2} (Khintchine's inequality) to $\tilde{\E} |I(f)|^4$ (i.e. to the randomness coming from $(f(p))_{x^{1-3\beta} < p \leq x^{1-\beta}}$ only) immediately implies that
$$ \E|I(f)|^4 \ll x^2 \E\left( \sum_{x^{1-3\beta} < p \leq x^{1-\beta}} \frac{1}{p} \Biggl| \int_{-\mathcal{T}}^{\mathcal{T}} \frac{1}{p^{iv}} F_{f}^{*}(1/2 + iv) G_{f}(1/2 + iv) \frac{x^{iv}}{1/2 + iv} dv \Biggr|^2 \right)^2 . $$
Note that, so far, we have used almost none of our assumptions or information about the parameters $\beta, \mathcal{T}, P$ and the functions $F_{f}^{*}(1/2 + iv),  G_{f}(1/2 + iv)$.

Again, rather than immediately expanding the square on the inside, we can produce a much cleaner argument with some preliminary work. This time, the work is of a harmonic analysis flavour. Slightly enlarging and smoothing the sum over $p$, we find
\begin{eqnarray}\label{eqafterkh}
&& \sum_{x^{1-3\beta} < p \leq x^{1-\beta}} \frac{1}{p} \Biggl| \int_{-\mathcal{T}}^{\mathcal{T}} \frac{1}{p^{iv}} F_{f}^{*}(1/2 + iv) G_{f}(1/2 + iv) \frac{x^{iv}}{1/2 + iv} dv \Biggr|^2 \\
& \ll & \frac{1}{\log x} \sum_{\substack{x^{1-4\beta} \\ < p \leq x}} \frac{\log p}{p} \left(1 - \frac{|\log p - (1-2\beta)\log x|}{2\beta\log x} \right) \Biggl| \int_{-\mathcal{T}}^{\mathcal{T}} \frac{1}{p^{iv}} F_{f}^{*}(\frac{1}{2} + iv) G_{f}(\frac{1}{2} + iv) \frac{x^{iv}}{1/2 + iv} dv \Biggr|^2 . \nonumber
\end{eqnarray}
Applying partial summation to the Prime Number Theorem with classical error term, for any $h \in \R$ we see
\begin{eqnarray}
&& \frac{1}{\log x} \sum_{x^{1-4\beta} < p \leq x} \frac{\log p}{p^{1+ih}} \left(1 - \frac{|\log p - (1-2\beta)\log x|}{2\beta\log x} \right) \nonumber \\
& = & \frac{1}{\log x} \int_{x^{1-4\beta}}^{x} \frac{1}{t^{1+ih}} \left(1 - \frac{|\log t - (1-2\beta)\log x|}{2\beta\log x} \right) dt + O((1+|h|) e^{-10c\sqrt{\log x}}) , \nonumber
\end{eqnarray}
say. Substituting $u=\log t$, the integral here is
\begin{eqnarray}
& = & \frac{1}{\log x} \int_{(1-4\beta)\log x}^{\log x} e^{-ihu} \left(1 - \frac{|u - (1-2\beta)\log x|}{2\beta\log x} \right) du \nonumber \\
& = & \frac{e^{-ih(1-2\beta)\log x}}{\log x} \int_{-2\beta\log x}^{2\beta\log x} e^{-ihw} \left(1 - \frac{|w|}{2\beta\log x} \right) dw , \nonumber
\end{eqnarray}
and Harmonic Analysis Result \ref{harman1} shows this is all $\ll \min\{\beta, \frac{1}{h^2 \beta\log^{2}x} \}$. If $|h| \leq e^{c\sqrt{\log x}}$, then the error term $O((1+|h|) e^{-10c\sqrt{\log x}})$ can be subsumed into this bound as well. Now opening the square on the second line of \eqref{eqafterkh}, and applying our bound with $|h| = |v_1 - v_2| \leq 2\mathcal{T}$, we deduce that
\begin{eqnarray}
&& \sum_{x^{1-3\beta} < p \leq x^{1-\beta}} \frac{1}{p} \Biggl| \int_{-\mathcal{T}}^{\mathcal{T}} \frac{1}{p^{iv}} F_{f}^{*}(1/2 + iv) G_{f}(1/2 + iv) \frac{x^{iv}}{1/2 + iv} dv \Biggr|^2 \nonumber \\
& \ll & \int_{-\mathcal{T}}^{\mathcal{T}} \int_{-\mathcal{T}}^{\mathcal{T}} \frac{|F_{f}^{*}(\frac{1}{2} + iv_1) G_{f}(\frac{1}{2} + iv_1)|}{|1/2 + iv_1|} \frac{|F_{f}^{*}(\frac{1}{2} + iv_2) G_{f}(\frac{1}{2} + iv_2)|}{|1/2 + iv_2|} \min\{\beta, \frac{1}{|v_1 - v_2|^2 \beta\log^{2}x} \} dv_1 dv_2 . \nonumber
\end{eqnarray}
Since $\frac{|F_{f}^{*}(\frac{1}{2} + iv_1) G_{f}(\frac{1}{2} + iv_1)|}{|1/2 + iv_1|} \frac{|F_{f}^{*}(\frac{1}{2} + iv_2) G_{f}(\frac{1}{2} + iv_2)|}{|1/2 + iv_2|} \leq  \frac{|F_{f}^{*}(\frac{1}{2} + iv_1)|^2 |G_{f}(\frac{1}{2} + iv_1)|^2}{|1/2 + iv_1|^2} + \frac{|F_{f}^{*}(\frac{1}{2} + iv_2)|^2 |G_{f}(\frac{1}{2} + iv_2)|^2}{|1/2 + iv_2|^2}$, and also using the resulting symmetry of $v_1$ and $v_2$, we find this is all
\begin{eqnarray}
& \ll & \int_{-\mathcal{T}}^{\mathcal{T}} \int_{-\mathcal{T}}^{\mathcal{T}} \frac{|F_{f}^{*}(\frac{1}{2} + iv_1)|^2 |G_{f}(\frac{1}{2} + iv_1)|^2}{|1/2 + iv_1|^2} \min\{\beta, \frac{1}{|v_1 - v_2|^2 \beta\log^{2}x} \} dv_1 dv_2 \nonumber \\
& \ll & \frac{1}{\log x} \int_{-\mathcal{T}}^{\mathcal{T}} \frac{|F_{f}^{*}(\frac{1}{2} + iv)|^2 |G_{f}(\frac{1}{2} + iv)|^2}{|1/2 + iv|^2} dv . \nonumber
\end{eqnarray}
Notice the role of the quadratic decay $|v_1 - v_2|^2$ here, which we introduced via the Fej\'er kernel. If one only had the analogous bound with linear decay, the overall estimate would be worsened by (at least) a factor of $\log\log x$, preventing us from deducing a sharp moment estimate. 

In summary, we have shown so far that
\begin{equation}\label{eqfourthafterred}
\E|I(f)|^4 \ll x^2 \E\left( \frac{1}{\log x} \int_{-\mathcal{T}}^{\mathcal{T}} \frac{|F_{f}^{*}(\frac{1}{2} + iv)|^2 |G_{f}(\frac{1}{2} + iv)|^2}{|1/2 + iv|^2} dv \right)^2 ,
\end{equation}
where we recall that
$$ F_{f}^{*}(s) := \Biggl(\sum_{\substack{n \leq x^{\beta}, \\ n \; \text{is} \; P \; \text{smooth}}} \frac{f(n)}{n^s} \Biggr) \cdot \Biggl(\sum_{\substack{n \leq x^{2\beta}, \\ n \; \text{is} \; P \; \text{rough}}} \frac{f(n)}{n^s} \Biggr) . $$

It will be convenient now, before working further with $F_{f}^{*}(1/2 + iv)$, to replace $G_{f}(1/2 + iv)$ with the more agreeable weight $W_{f}(1/2 + iv)$. To do this, we may note that $|G_{f}(1/2 + iv)|^2 \ll |W_{f}(1/2 + iv)|^2 + |G_{f}(1/2 + iv) - W_{f}(1/2 + iv)|^2$, and (applying the Cauchy--Schwarz inequality to the integral, and using the bound $\int_{-\mathcal{T}}^{\mathcal{T}} \frac{dv}{|1/2 + iv|^2} \ll 1$) that
\begin{eqnarray}
&& \frac{x^2}{\log^{2}x} \E\left( \int_{-\mathcal{T}}^{\mathcal{T}} \frac{|F_{f}^{*}(\frac{1}{2} + iv)|^2 |G_{f}(\frac{1}{2} + iv) - W_{f}(\frac{1}{2} + iv)|^2}{|1/2 + iv|^2} dv \right)^2 \nonumber \\
& \ll & \frac{x^2}{\log^{2}x} \int_{-\mathcal{T}}^{\mathcal{T}} \frac{\E|F_{f}^{*}(\frac{1}{2} + iv)|^4 |G_{f}(\frac{1}{2} + iv) - W_{f}(\frac{1}{2} + iv)|^4}{|1/2 + iv|^2} dv . \nonumber
\end{eqnarray}
This is extremely similar to the quantities we dealt with in section \ref{subsecgtow}. Using the independence of $\sum_{\substack{n \leq x^{2\beta}, \\ n \; \text{is} \; P \; \text{rough}}} \frac{f(n)}{n^{1/2+iv}}$ from $\sum_{\substack{n \leq x^{\beta}, \\ n \; \text{is} \; P \; \text{smooth}}} \frac{f(n)}{n^{1/2+iv}}, G_{f}(1/2 + iv), W_{f}(1/2 + iv)$, we find $\E|F_{f}^{*}(1/2 + iv)|^4 |G_{f}(1/2 + iv) - W_{f}(1/2 + iv)|^4$ is
$$ = \E|\sum_{\substack{n \leq x^{2\beta}, \\ n \; \text{is} \; P \; \text{rough}}} \frac{f(n)}{n^{1/2+iv}}|^4 \cdot \E|\sum_{\substack{n \leq x^{\beta}, \\ n \; \text{is} \; P \; \text{smooth}}} \frac{f(n)}{n^{1/2+iv}}|^4 |G_{f}(\frac{1}{2} + iv) - W_{f}(\frac{1}{2} + iv)|^4 , $$
and then the Cauchy--Schwarz inequality and Probability Result \ref{probres1} imply this is
\begin{eqnarray}
& \leq & \E|\sum_{\substack{n \leq x^{2\beta}, \\ n \; \text{is} \; P \; \text{rough}}} \frac{f(n)}{n^{1/2+iv}}|^4 \cdot \sqrt{\E|\sum_{\substack{n \leq x^{\beta}, \\ n \; \text{is} \; P \; \text{smooth}}} \frac{f(n)}{n^{1/2+iv}}|^8 } \sqrt{\E|G_{f}(\frac{1}{2} + iv) - W_{f}(\frac{1}{2} + iv)|^8} \nonumber \\
& \ll & \Biggl( \sum_{\substack{n \leq x^{2\beta}, \\ n \; \text{is} \; P \; \text{rough}}} \frac{d(n)}{n} \Biggr)^2 \Biggl( \sum_{\substack{n \leq x^{\beta}, \\ n \; \text{is} \; P \; \text{smooth}}} \frac{d_{4}(n)}{n} \Biggr)^2 \sqrt{\E|G_{f}(\frac{1}{2} + iv) - W_{f}(\frac{1}{2} + iv)|^8 } . \nonumber
\end{eqnarray}
Here we have $\sum_{\substack{n \leq x^{2\beta}, \\ n \; \text{is} \; P \; \text{rough}}} \frac{d(n)}{n} \leq \prod_{P < p \leq x^{2\beta}} (1 - \frac{1}{p})^{-2} \ll (\frac{\log(x^{\beta})}{\log P})^2$, and $\sum_{\substack{n \leq x^{\beta}, \\ n \; \text{is} \; P \; \text{smooth}}} \frac{d_{4}(n)}{n} \leq \prod_{p \leq P} (1 - \frac{1}{p})^{-4} \ll \log^{4}P$. Meanwhile, identical calculations as in section \ref{subsecgtow} yield that $\sqrt{\E|G_{f}(\frac{1}{2} + iv) - W_{f}(\frac{1}{2} + iv)|^8} \ll \frac{1}{(\log\log P)^{(\log\log P)^{10}}} \ll \frac{1}{\log^{1004}P}$, say. In conjunction with \eqref{eqfourthafterred}, we then arrive at the bound
$$ \E|I(f)|^4 \ll x^2 \E\left( \frac{1}{\log x} \int_{-\mathcal{T}}^{\mathcal{T}} \frac{|F_{f}^{*}(\frac{1}{2} + iv)|^2 |W_{f}(\frac{1}{2} + iv)|^2}{|1/2 + iv|^2} dv \right)^2 + \frac{\beta^4 x^2 \log^{2}x}{\log^{1000}P} , $$
where the second term is clearly more than good enough.

As our final preparatory step, we need to replace the smooth sum $\sum_{\substack{n \leq x^{\beta}, \\ n \; \text{is} \; P \; \text{smooth}}} \frac{f(n)}{n^s}$ inside $F_{f}^{*}(s)$ by the Euler product $F_{P}(s) = \prod_{p \leq P} (1 - \frac{f(p)}{p^s})^{-1} = \sum_{\substack{n=1, \\ P \; \text{smooth}}}^{\infty} \frac{f(n)}{n^{s}}$, similarly as in the proof of Key Proposition \ref{keyprop3}. The error that we make in doing this (and then applying the Cauchy--Schwarz inequality to the resulting integral over $v$, as above) will be
\begin{equation}\label{eqgotoFP}
\ll \frac{x^2}{\log^{2}x} \int_{-\mathcal{T}}^{\mathcal{T}} \frac{\E|\sum_{\substack{n \leq x^{2\beta}, \\ n \; \text{is} \; P \; \text{rough}}} \frac{f(n)}{n^{1/2 + iv}}|^4 |\sum_{\substack{n > x^{\beta}, \\ n \; \text{is} \; P \; \text{smooth}}} \frac{f(n)}{n^{1/2 + iv}}|^4 |W_{f}(\frac{1}{2} + iv)|^4}{|1/2 + iv|^2} dv .
\end{equation}
Using the uniform boundedness of $W_f$ (which would not be available for $G_f$); the independence of $f(n)$ on $P$-smooth and $P$-rough numbers; and then Probability Result \ref{probres1}; we find the numerator is
$$ \ll \E|\sum_{\substack{n \leq x^{2\beta}, \\ n \; \text{is} \; P \; \text{rough}}} \frac{f(n)}{n^{1/2 + iv}}|^4 \cdot \E|\sum_{\substack{n > x^{\beta}, \\ n \; \text{is} \; P \; \text{smooth}}} \frac{f(n)}{n^{1/2 + iv}}|^4 \leq \Biggl( \sum_{\substack{n \leq x^{2\beta}, \\ n \; \text{is} \; P \; \text{rough}}} \frac{d(n)}{n} \Biggr)^2 \Biggl( \sum_{\substack{n > x^{\beta}, \\ n \; \text{is} \; P \; \text{smooth}}} \frac{d(n)}{n} \Biggr)^2 . $$
As above, we have $\sum_{\substack{n \leq x^{2\beta}, \\ n \; \text{is} \; P \; \text{rough}}} \frac{d(n)}{n} \ll (\frac{\log(x^{2\beta})}{\log P})^2$. And using Rankin's trick, we get
$$ \sum_{\substack{n > x^{\beta}, \\ P \; \text{smooth}}} \frac{d(n)}{n} \leq x^{-\beta/\log P} \sum_{\substack{n = 1, \\ P \; \text{smooth}}}^{\infty} \frac{d(n)}{n^{1-1/\log P}} = x^{-\beta/\log P} \prod_{p \leq P}(1 - \frac{1}{p^{1-1/\log P}})^{-2} \ll x^{-\beta/\log P} \log^{2}P . $$
So we conclude that \eqref{eqgotoFP} is $\ll \frac{\beta^4 x^2 \log^{2}x}{x^{2\beta/\log P}}$.

\subsection{Exploitation of the barrier}\label{subsecexploitbarr}
We have now shown that
\begin{eqnarray}
\E|I(f)|^4 & \ll & x^2 \E\Biggl( \frac{1}{\log x} \int_{-\mathcal{T}}^{\mathcal{T}} \frac{|\sum_{\substack{n \leq x^{2\beta}, \\ n \; \text{is} \; P \; \text{rough}}} \frac{f(n)}{n^{1/2 + iv}}|^2 |F_{P}(\frac{1}{2} + iv)|^2 |W_{f}(\frac{1}{2} + iv)|^2}{|1/2 + iv|^2} dv \Biggr)^2 + \nonumber \\
&& + \frac{\beta^4 x^2 \log^{2}x}{\log^{1000}P} + \frac{\beta^4 x^2 \log^{2}x}{x^{2\beta/\log P}} . \nonumber
\end{eqnarray}

Before dealing with the expectation here in earnest, we do a little more tidying up. Breaking the range of integration over $v$ into sub-intervals, and then applying the Cauchy--Schwarz inequality, we can upper bound the first line of the previous display by
\begin{eqnarray}
& \ll & \frac{x^2}{\log^{2}x} \E\Biggl( \sum_{V=-\infty}^{\infty} \frac{1}{1+V^2} \int_{V-1/2}^{V+1/2} |\sum_{\substack{n \leq x^{2\beta}, \\ n \; \text{is} \; P \; \text{rough}}} \frac{f(n)}{n^{1/2 + iv}}|^2 |F_{P}(\frac{1}{2} + iv)|^2 |W_{f}(\frac{1}{2} + iv)|^2 dv \Biggr)^2 \nonumber \\
& \ll & \frac{x^2}{\log^{2}x} \sum_{V=-\infty}^{\infty} \frac{1}{1+V^2} \E\Biggl( \int_{V-1/2}^{V+1/2} |\sum_{\substack{n \leq x^{2\beta}, \\ n \; \text{is} \; P \; \text{rough}}} \frac{f(n)}{n^{1/2 + iv}}|^2 |F_{P}(\frac{1}{2} + iv)|^2 |W_{f}(\frac{1}{2} + iv)|^2 dv \Biggr)^2 . \nonumber
\end{eqnarray}
Translation invariance in law implies that the expectation on the inside is independent of $V$, and so it will suffice to bound it when $V=0$, say. Recalling \eqref{eqWsharper}, and then applying Probability Result \ref{probres1} as in many previous calculations, we see that contribution is
\begin{eqnarray}\label{eqlikeLM}
& \ll & \frac{x^2}{\log^{2}x} \E\Biggl( \int_{-1/2}^{1/2} |\sum_{\substack{n \leq x^{2\beta}, \\ n \; \text{is} \; P \; \text{rough}}} \frac{f(n)}{n^{1/2 + iv}}|^2 |F_{P}(\frac{1}{2} + iv)|^2 (\textbf{1}_{U_{f}(v)} + \frac{1}{\log^{2000}P}) dv \Biggr)^2 \nonumber \\
& \ll & \frac{x^2}{\log^{2}x} \E\Biggl( \int_{-1/2}^{1/2} |\sum_{\substack{n \leq x^{2\beta}, \\ n \; \text{is} \; P \; \text{rough}}} \frac{f(n)}{n^{1/2 + iv}}|^2 |F_{P}(\frac{1}{2} + iv)|^2 \textbf{1}_{U_{f}(v)} dv \Biggr)^2 + \\
&& + \frac{x^2}{\log^{2}x} \frac{1}{\log^{4000}P} \Biggl( \sum_{\substack{n \leq x^{2\beta}, \\ n \; \text{is} \; P \; \text{rough}}} \frac{d(n)}{n} \Biggr)^2 \Biggl( \sum_{\substack{n = 1, \\ n \; \text{is} \; P \; \text{smooth}}}^{\infty} \frac{d(n)}{n} \Biggr)^2 . \nonumber
\end{eqnarray}
The second line is $\ll \frac{\beta^4 x^2 \log^{2}x}{\log^{4000}P}$, which is more than acceptable.

\vspace{12pt}
For our main argument, we can essentially follow the proof of Key Proposition 5 of Harper~\cite{harperrmflow}, with some extra care to handle the contribution from the rough sum $\sum_{\substack{n \leq x^{2\beta}, \\ n \; \text{is} \; P \; \text{rough}}} \frac{f(n)}{n^{1/2 + iv}}$ (which does not appear in the purely probabilistic setting there). For ease of checking, and since this is the real heart of the argument and place where the choice of the barrier $U_{f}$ (and thus of the weight $W$) becomes crucial, we shall give full details.

Expanding out, we find the expectation on the first line of \eqref{eqlikeLM} is
$$ = \int_{-1/2}^{1/2} \int_{-1/2}^{1/2} \E \textbf{1}_{U_{f}(v)} |\sum_{\substack{n \leq x^{2\beta}, \\ P \; \text{rough}}} \frac{f(n)}{n^{1/2 + iv}}|^2 |F_{P}(\frac{1}{2} + iv)|^2 \textbf{1}_{U_{f}(w)} |\sum_{\substack{n \leq x^{2\beta}, \\ P \; \text{rough}}} \frac{f(n)}{n^{1/2 + iw}}|^2 |F_{P}(\frac{1}{2} + iw)|^2 dv dw . $$
Again, we can use translation invariance in law to simplify this by shifting $v$ to 0, and replacing $w$ by $w-v$. This yields an upper bound
\begin{equation}\label{upperaftertr}
\leq \int_{-1}^{1} \E \textbf{1}_{U_{f}(0)} |\sum_{\substack{n \leq x^{2\beta}, \\ P \; \text{rough}}} \frac{f(n)}{n^{1/2}}|^2 |F_{P}(\frac{1}{2})|^2 \textbf{1}_{U_{f}(w)} |\sum_{\substack{n \leq x^{2\beta}, \\ P \; \text{rough}}} \frac{f(n)}{n^{1/2 + iw}}|^2 |F_{P}(\frac{1}{2} + iw)|^2 dw .
\end{equation}

Now if $U_{f}(w)$ (as defined at the beginning of section \ref{secpasstormf}, with $f$ replacing $\chi$) occurs, and we temporarily set $J = J(w) := \max\{\lfloor \log(|w|\log P) \rfloor, 0\}$, then we have
\begin{eqnarray}\label{lwoccurs}
|F_{P^{e^{-J}}}(\frac{1}{2} + iw)|^{2} & \ll & \Biggl( \frac{\log P}{e^{J} (\log\log P - J)^{50}} e^{\min\{\sqrt{\log\log P}, \frac{1}{1-q}\}} \Biggr)^2 \\
& \ll & \Biggl( \min\{\frac{\log P}{(\log\log P)^{50}} , \frac{1}{|w|\log^{50}(2/|w|)}\} e^{\min\{\sqrt{\log\log P}, \frac{1}{1-q}\}} \Biggr)^2 . \nonumber
\end{eqnarray}
Note that if $J(w) > \log\log P - 1$, then the condition $U_{f}(w)$ is not applicable but $P^{e^{-J}} \leq e^{e}$, so \eqref{lwoccurs} certainly still holds (indeed the left hand side is uniformly bounded). 

It particular, when $|w| \leq 1/(\log P)^{1/3}$ (say) we get a bound
$$ |F_{P^{e^{-J}}}(\frac{1}{2} + iw)|^{2} \ll \frac{1}{(\log\log P)^{100}} \Biggl( \min\{\log P , \frac{1}{|w|}\} e^{\min\{\sqrt{\log\log P}, \frac{1}{1-q}\}} \Biggr)^2 . $$
Here the saving $1/(\log\log P)^{100}$ will be sufficient that, on this range of $w$, we have no further need of the barrier events $\textbf{1}_{U_{f}(0)}, \textbf{1}_{U_{f}(w)}$. Discarding these, and using the independence of $f$ on different primes, the expectation inside our integral over $w$ in \eqref{upperaftertr} is
\begin{eqnarray}
& \ll & \frac{1}{(\log\log P)^{100}} \Biggl( \min\{\log P , \frac{1}{|w|}\} e^{\min\{\sqrt{\log\log P}, \frac{1}{1-q}\}} \Biggr)^2 \cdot \E |\sum_{\substack{n \leq x^{2\beta}, \\ P \; \text{rough}}} \frac{f(n)}{n^{1/2}}|^2 |\sum_{\substack{n \leq x^{2\beta}, \\ P \; \text{rough}}} \frac{f(n)}{n^{1/2 + iw}}|^2 \cdot \nonumber \\
&& \cdot \E |F_{P^{e^{-J}}}(\frac{1}{2})|^2 \cdot \E \prod_{P^{e^{-J}} < p \leq P} |1 - \frac{f(p)}{p^{1/2}}|^{-2} \cdot \prod_{P^{e^{-J}} < p \leq P} |1 - \frac{f(p)}{p^{1/2 + iw}}|^{-2} . \nonumber
\end{eqnarray}
As we have seen before, the Cauchy--Schwarz inequality and Probability Result \ref{probres1} imply that $\E |\sum_{\substack{n \leq x^{2\beta}, \\ n \; \text{is} \; P \; \text{rough}}} \frac{f(n)}{n^{1/2}}|^2 |\sum_{\substack{n \leq x^{2\beta}, \\ n \; \text{is} \; P \; \text{rough}}} \frac{f(n)}{n^{1/2 + iw}}|^2 \leq ( \sum_{\substack{n \leq x^{2\beta}, \\ n \; \text{is} \; P \; \text{rough}}} \frac{d(n)}{n} )^2 \leq \prod_{P < p \leq x^{2\beta}} (1 - \frac{1}{p})^{-4} \ll (\frac{\log(x^{\beta})}{\log P})^4$. We also have $\E |F_{P^{e^{-J}}}(\frac{1}{2})|^2 \ll \frac{\log P}{e^J} \ll \min\{\log P , \frac{1}{|w|}\}$, and Probability Result \ref{probres5} implies (since $|w| \geq \frac{1}{\log(P^{e^{-J}})}$, except when $J=0$ and the products are empty anyway) that $\E \prod_{P^{e^{-J}} < p \leq P} |1 - \frac{f(p)}{p^{1/2}}|^{-2} \cdot \prod_{P^{e^{-J}} < p \leq P} |1 - \frac{f(p)}{p^{1/2 + iw}}|^{-2} \ll e^{2J} \ll \max\{(|w|\log P)^2 , 1\}$. So the total contribution to \eqref{upperaftertr} from $|w| \leq 1/(\log P)^{1/3}$ is
\begin{eqnarray}
& \ll & \frac{e^{2\min\{\sqrt{\log\log P}, \frac{1}{1-q}\}}}{(\log\log P)^{100}} \log^{2}P (\frac{\log(x^{\beta})}{\log P})^4 \int_{-1/(\log P)^{1/3}}^{1/(\log P)^{1/3}} \min\{\log P , \frac{1}{|w|}\} dw \nonumber \\
& \ll & \frac{e^{2\min\{\sqrt{\log\log P}, \frac{1}{1-q}\}}}{(\log\log P)^{99}} \log^{2}(x^{\beta}) (\frac{\log(x^{\beta})}{\log P})^2 , \nonumber
\end{eqnarray}
which is acceptable.

On the remaining range $1/(\log P)^{1/3} < |w| \leq 1$, Probability Result \ref{probres3} importantly implies the stronger (sharp) bound $\E |\sum_{\substack{n \leq x^{2\beta}, \\ n \; \text{is} \; P \; \text{rough}}} \frac{f(n)}{n^{1/2}}|^2 |\sum_{\substack{n \leq x^{2\beta}, \\ n \; \text{is} \; P \; \text{rough}}} \frac{f(n)}{n^{1/2 + iw}}|^2 \ll (\frac{\log(x^{\beta})}{\log P})^2$ for the rough number contribution. Invoking this, along with \eqref{lwoccurs} and the independence of $f$ on different primes, we see the contribution from $1/(\log P)^{1/3} < |w| \leq 1$ is
\begin{eqnarray}
& \ll & e^{2\min\{\sqrt{\log\log P}, \frac{1}{1-q}\}} (\frac{\log(x^{\beta})}{\log P})^2 \cdot \nonumber \\
&& \cdot \int_{\frac{1}{(\log P)^{1/3}} < |w| \leq 1} \frac{1}{|w|^2 \log^{100}(2/|w|)} \E \textbf{1}_{U_{f}(0)} \prod_{p \leq P} |1 - \frac{f(p)}{p^{1/2}}|^{-2} \cdot \textbf{1}_{U_{f}(w)} \prod_{\substack{P^{e^{-J}} \\ < p \leq P}} |1 - \frac{f(p)}{p^{1/2 + iw}}|^{-2} dw . \nonumber
\end{eqnarray}

This time, we need to make further use of the indicator functions $\textbf{1}_{U_{f}(0)}, \textbf{1}_{U_{f}(w)}$ to obtain a sharp bound, because our saving $1/\log^{100}(2/|w|)$ from \eqref{lwoccurs} is no longer so great (especially when $|w| \asymp 1$, which will ultimately contribute the ``main term''). We temporarily set $K = K(w) := \lfloor \log(|w|^2\log P) - 2B \rfloor$, where $B$ is the large absolute constant from Probability Result \ref{probres5}. Note that this choice satisfies $K(w) \asymp \log\log P$, and $\frac{\log P}{e^K} \asymp \frac{1}{|w|^2}$, and $K(w) \leq J(w)$ on our range $1/(\log P)^{1/3} < |w| \leq 1$. If the event $U_{f}(0)$ occurs then, comparing the definition of $U_{f}(0)$ with $j = K$ and with general $j$, we must in particular have
\begin{eqnarray}
\Biggl(\frac{\log P}{e^{j}} \Biggr)^{-B} e^{-2\min\{\sqrt{\log\log P}, \frac{1}{1-q}\}} |w|^2 \log^{50}(\frac{2}{|w|}) & \ll & \prod_{P^{e^{-K}} < p \leq P^{e^{-j}}} |1 - \frac{f(p)}{p^{1/2}}|^{-1} \nonumber \\
& \ll & \frac{\log P}{e^{j}} e^{- 50\log\log(\frac{\log P}{e^{j}})} \frac{e^{2\min\{\sqrt{\log\log P}, 1/(1-q)\}}}{|w|^{2B}} \nonumber
\end{eqnarray}
for all $0 \leq j \leq K-1$. Similarly, if the event $U_{f}(w)$ occurs then we must have the same bounds with $\prod_{P^{e^{-K}} < p \leq P^{e^{-j}}} |1 - \frac{f(p)}{p^{1/2}}|^{-1}$ replaced by $\prod_{P^{e^{-K}} < p \leq P^{e^{-j}}} |1 - \frac{f(p)}{p^{1/2 + iw}}|^{-1}$. Let us write $R_{f}(w)$ for the event that one has these bounds for both $\prod_{P^{e^{-K}} < p \leq P^{e^{-j}}} |1 - \frac{f(p)}{p^{1/2}}|^{-1}$ and $\prod_{P^{e^{-K}} < p \leq P^{e^{-j}}} |1 - \frac{f(p)}{p^{1/2 + iw}}|^{-1}$ simultaneously. Then we can upper bound the expectation remaining in our integral by
$$ \E \textbf{1}_{R_{f}(w)} \prod_{p \leq P} |1 - \frac{f(p)}{p^{1/2}}|^{-2} \cdot \prod_{P^{e^{-J}} < p \leq P} |1 - \frac{f(p)}{p^{1/2 + iw}}|^{-2} . $$

Since $f(p)$ is independent on distinct primes, and the event $R_{f}(w)$ only depends on the values of $f$ on primes $> P^{e^{-K}} \geq P^{e^{-J}}$, this expectation factors as
\begin{eqnarray}
&& \E |F_{P^{e^{-J}}}(\frac{1}{2})|^2 \cdot \E \prod_{P^{e^{-J}} < p \leq P^{e^{-K}}} |1 - \frac{f(p)}{p^{1/2}}|^{-2} \cdot \prod_{P^{e^{-J}} < p \leq P^{e^{-K}}} |1 - \frac{f(p)}{p^{1/2 + iw}}|^{-2} \cdot \nonumber \\
&& \cdot \E \textbf{1}_{R_{f}(w)} \prod_{P^{e^{-K}} < p \leq P} |1 - \frac{f(p)}{p^{1/2}}|^{-2} \cdot \prod_{P^{e^{-K}} < p \leq P} |1 - \frac{f(p)}{p^{1/2 + iw}}|^{-2} . \nonumber
\end{eqnarray}
We have $\E |F_{P^{e^{-J}}}(\frac{1}{2})|^2 \ll \frac{\log P}{e^J}$, and the first part of Probability Result \ref{probres5} implies (noting $|w| \geq \frac{1}{\log(P^{e^{-J}})}$) that the expectation of the products over $P^{e^{-J}} < p \leq P^{e^{-K}}$ is $\ll e^{2(J-K)}$. Finally, the second part of Probability Result \ref{probres5} (which is applicable with $a = 2\min\{\sqrt{\log\log P}, \frac{1}{1-q}\} + O(\log(2/|w|))$ and $h(y) = -50\log y$, since $\log(P^{e^{-K}}) \geq \frac{e^{2B}}{|w|^2}$) implies that
\begin{eqnarray}
&& \E \textbf{1}_{R_{f}(w)} \prod_{P^{e^{-K}} < p \leq P} |1 - \frac{f(p)}{p^{1/2}}|^{-2} \cdot \prod_{P^{e^{-K}} < p \leq P} |1 - \frac{f(p)}{p^{1/2 + iw}}|^{-2} \nonumber \\
& \ll & \min\left\{1, \frac{\min\{\sqrt{\log\log P}, 1/(1-q)\} + \log(2/|w|)}{\sqrt{K}} \right\}^2 e^{2K} \nonumber \\
& \ll & \log^{2}(2/|w|) \min\left\{1, \frac{1}{(1-q)\sqrt{\log\log P}} \right\}^2 e^{2K} , \nonumber
\end{eqnarray}
on recalling that $\sqrt{K} \asymp \sqrt{\log\log P}$.

Thus the overall contribution to \eqref{upperaftertr} from $1/(\log P)^{1/3} < |w| \leq 1$ is
\begin{eqnarray}
& \ll & e^{2\min\{\sqrt{\log\log P}, \frac{1}{1-q}\}} \min\left\{1, \frac{1}{(1-q)\sqrt{\log\log P}} \right\}^2 (\frac{\log(x^{\beta})}{\log P})^2 \log P \cdot \nonumber \\
&& \cdot \int_{\frac{1}{(\log P)^{1/3}} < |w| \leq 1} \frac{1}{|w|^2 \log^{98}(2/|w|)} e^{J(w)} dw . \nonumber
\end{eqnarray}
Recalling $e^{J(w)} \asymp |w|\log P$ here, we find the integral is $\ll \log P \int_{\frac{1}{(\log P)^{1/3}} < |w| \leq 1} \frac{dw}{|w| \log^{98}(2/|w|)} \ll \log P$ and the total contribution is $\ll e^{2\min\{\sqrt{\log\log P}, \frac{1}{1-q}\}} \min\left\{1, \frac{1}{(1-q)\sqrt{\log\log P}} \right\}^2 (\beta \log x)^2$, which is also acceptable when plugged into \eqref{eqlikeLM} and multiplied by $\frac{x^2}{\log^{2}x}$. We emphasise the quite delicate nature of this calculation, and in particular the crucial role of the term $ \log^{98}(2/|w|)$ in the denominator in making the final integral uniformly bounded. The presence of this flows (via \eqref{lwoccurs}) from the term $- 50\log\log(\frac{\log P}{e^{j}})$ in our initial barrier construction in section \ref{secpasstormf}, and is one of the key reasons for including such a term there (the other being to compensate when $|w| \leq 1/(\log P)^{1/3}$ for the extra powers of $(\frac{\log(x^{\beta})}{\log P})$ in our bounds, arising from the lack of rough number decorrelation on that range).
\qed

\subsection{The second moment case}
Proving the second moment part of Key Proposition \ref{keyprop4} will be a significantly simpler variant of the fourth moment proof.

Proceeding similarly as in section \ref{subsecfourthkh}, but now writing $\E|I(f)|^2 = \E\tilde{\E} |I(f)|^2$ and simply applying orthogonality (rather than Khintchine's inequality) to $\tilde{\E} |I(f)|^2$, we get
$$ \E|I(f)|^2 \ll x \E\left( \sum_{x^{1-3\beta} < p \leq x^{1-\beta}} \frac{1}{p} \Biggl| \int_{-\mathcal{T}}^{\mathcal{T}} \frac{1}{p^{iv}} F_{f}^{*}(1/2 + iv) G_{f}(1/2 + iv) \frac{x^{iv}}{1/2 + iv} dv \Biggr|^2 \right) . $$
Upper bounding the sum over primes as in \eqref{eqafterkh}, we ultimately deduce that
$$ \E|I(f)|^2 \ll \frac{x}{\log x} \E \int_{-\mathcal{T}}^{\mathcal{T}} \frac{|F_{f}^{*}(\frac{1}{2} + iv)|^2 |G_{f}(\frac{1}{2} + iv)|^2}{|1/2 + iv|^2} dv . $$

Since
$$ F_{f}^{*}(1/2 + iv) := \Biggl(\sum_{\substack{n \leq x^{\beta}, \\ n \; \text{is} \; P \; \text{smooth}}} \frac{f(n)}{n^{1/2 + iv}} \Biggr) \cdot \Biggl(\sum_{\substack{n \leq x^{2\beta}, \\ n \; \text{is} \; P \; \text{rough}}} \frac{f(n)}{n^{1/2 + iv}} \Biggr) , $$
with the $P$-rough sum independent of the $P$-smooth sum and of $G_{f}(1/2 + iv)$, we can pull out a factor of $\E|\sum_{\substack{n \leq x^{2\beta}, \\ n \; \text{is} \; P \; \text{rough}}} \frac{f(n)}{n^{1/2+iv}}|^2 = \sum_{\substack{n \leq x^{2\beta}, \\ n \; \text{is} \; P \; \text{rough}}} \frac{1}{n} \ll \frac{\log(x^{2\beta})}{\log P}$ and find that
$$ \E|I(f)|^2 \ll \frac{\beta x}{\log P} \E \int_{-\mathcal{T}}^{\mathcal{T}} \frac{|\sum_{\substack{n \leq x^{\beta}, \\ n \; \text{is} \; P \; \text{smooth}}} \frac{f(n)}{n^{1/2 + iv}}|^2 |G_{f}(\frac{1}{2} + iv)|^2}{|1/2 + iv|^2} dv . $$

Finally, similar calculations as in section \ref{subsecgtow} or section \ref{subsecfourthkh} show that we can replace $G_{f}(\frac{1}{2} + iv)$ with $W_{f}(\frac{1}{2} + iv)$, at the cost of an additive error term that is $\ll \frac{\beta x}{\log^{1000}P}$ (say); and then replace $\sum_{\substack{n \leq x^{\beta}, \\ n \; \text{is} \; P \; \text{smooth}}} \frac{f(n)}{n^{1/2 + iv}}$ with $F_{P}(\frac{1}{2} + iv)$, at the cost of an error term that is $\ll \frac{\beta x}{x^{\beta/\log P}}$. Recalling also \eqref{eqWsharper}, we then arrive at the bound
$$ \E|I(f)|^2 \ll \frac{\beta x}{\log P} \int_{-\mathcal{T}}^{\mathcal{T}} \frac{\E|F_{P}(1/2 + iv)|^2 (\textbf{1}_{U_{f}(v)} + 1/\log^{2000}P)}{|1/2 + iv|^2} dv + \frac{\beta x}{\log^{1000}P} + \frac{\beta x}{x^{\beta/\log P}} . $$
The term $\frac{\beta x}{\log^{1000}P}$ is more than good enough. Using the easy bound $\E|F_{P}(1/2 + iv)|^2 \ll \log P$, we see the total contribution from the $1/\log^{2000}P$ term is $\ll \frac{\beta x}{\log^{2000}P}$, which is also more than good enough. Then translation invariance in law implies that $\int_{-\mathcal{T}}^{\mathcal{T}} \frac{\E|F_{P}(1/2 + iv)|^2 \textbf{1}_{U_{f}(v)}}{|1/2 + iv|^2} dv = \int_{-\mathcal{T}}^{\mathcal{T}} \frac{\E|F_{P}(1/2)|^2 \textbf{1}_{U_{f}(0)}}{|1/2 + iv|^2} dv \ll \E|F_{P}(1/2)|^2 \textbf{1}_{U_{f}(0)}$, and Probability Result \ref{probres4} with $a = \min\{\sqrt{\log\log P}, \frac{1}{1-q}\} + A + 3$ and $h(y) = -50\log y$ shows this is $\ll \frac{\log P}{1 + (1-q)\sqrt{\log\log P}}$. So the total contribution from this term is $\ll \frac{\beta x}{\log P} \frac{\log P}{1 + (1-q)\sqrt{\log\log P}} \ll \frac{\beta x}{1 + (1-q)\sqrt{\log\log P}}$, as desired.
\qed

\section{Proofs of Theorems \ref{thmlb2} and \ref{thmlb3}}

\subsection{The weighted case: Theorem \ref{thmlb2}}
This theorem, where one has a multiplicative (but not necessarily totally multiplicative) twist $h(n)$ in the character sum, may be proved by reasonably easy adaptation of the arguments for Theorem \ref{thmlb1} (very similarly to the way one adapts the upper bound \eqref{eqbtscupper} to the weighted case).

In place of $I(\chi)$, we simply take as our proxy object the integral
$$ \frac{\sqrt{x}}{2\pi} \int_{-\mathcal{T}}^{\mathcal{T}} \Biggl(\sum_{x^{1-3\beta} < p \leq x^{1-\beta}} \frac{h(p) \chi(p)}{p^{1/2+iv}} \Biggr) F_{\chi, h}^{*}(1/2 + iv) G_{\chi, h}(1/2 + iv) \frac{x^{iv}}{1/2 + iv} dv . $$
Here we let
$$ F_{\chi, h}^{*}(s) := \Biggl(\sum_{\substack{n \leq x^{\beta}, \\ n \; \text{is} \; P \; \text{smooth}}} \frac{h(n) \chi(n)}{n^s} \Biggr) \cdot \Biggl(\sum_{\substack{n \leq x^{2\beta}, \\ n \; \text{is} \; P \; \text{rough}}} \frac{h(n) \chi(n)}{n^s} \Biggr) , $$
defining $G_{\chi, h}(s)$ (and the corresponding $W_{\chi, h}(s)$) analogously to $G_{\chi}(s)$ and $W_{\chi}(s)$, but with $\Re \sum_{p \leq P^{e^{-j}}} \sum_{k=1,2} \frac{\chi(p^k)}{kp^{ks}}$ replaced by $\Re \sum_{5 \leq p \leq P^{e^{-j}}} (\frac{h(p) \chi(p)}{p^{s}} + \frac{(h(p^2) - (1/2)h(p)^2)\chi(p)^2}{p^{2s}})$. Each summand here corresponds to the first and second order terms in the Taylor expansion of $\log|1 + \frac{h(p) \chi(p)}{p^{1/2+it}} + \sum_{k=2}^{\infty} \frac{h(p^k) \chi(p)^k}{p^{k(1/2+it)}}|$, the logarithm of the $h(n)$-weighted Euler factor corresponding to $p$. Note that $|\frac{h(p) \chi(p)}{p^{1/2+it}} + \sum_{k=2}^{\infty} \frac{h(p^k) \chi(p)^k}{p^{k(1/2+it)}}| \leq \sum_{k=1}^{\infty} \frac{1}{p^{k/2}} = \frac{1}{\sqrt{p} - 1}$. Provided that $p \geq 5$, this is all $\leq \frac{1}{\sqrt{5} - 1} < 1$, so we can legitimately apply Taylor expansion. This is not the case for the primes 2 and 3, but for those we still have an upper bound $\ll 1$ for the Euler factors, as well as (on the random side) $\E|1 + \frac{h(p) f(p)}{p^{1/2+it}} + \sum_{k=2}^{\infty} \frac{h(p^k) f(p)^k}{p^{k(1/2+it)}}|^2 = 1 + \frac{|h(p)|^2}{p} + \sum_{k=2}^{\infty} \frac{|h(p^k)|^2}{p^{k}} \geq 1 + \frac{1}{p}$.

With this choice of proxy, and with $|\sum_{n \leq x} \chi(n)|$ replaced by $|\sum_{n \leq x} h(n) \chi(n)|$, the same bounds as in Key Propositions \ref{keyprop1} and \ref{keyprop2} may be shown to hold uniformly for all multiplicative functions $h(n)$ as in Theorem \ref{thmlb2}. Indeed, one may pass to averages of random multiplicative functions exactly as in section \ref{secpasstormf}. Since $|h(n)| \leq 1$ for all $n$, the handling of most error terms in the calculations is identical to our previous arguments, and since $|h(p)|=1$ on primes we arrive at exactly the same prime number sums $\sum_{x^{1-3\beta} < p \leq x^{1-\beta}} \frac{1}{p^{1+i(t-v)}}, \sum_{x^{1-3\beta} < p \leq x^{1-\beta}} \frac{1}{p}$ as well. When handling the $P$-rough sum inside $F_{\chi, h}^{*}(s)$, in the analogue of Key Propositions \ref{keyprop1} and \ref{keyprop3} one ends up working with $\sum_{\substack{n \leq x^{2\beta}, \\ n \; \text{is} \; P \; \text{rough}}} \frac{|h(n)|^2}{n^{1+iw}}$ in place of $\sum_{\substack{n \leq x^{2\beta}, \\ n \; \text{is} \; P \; \text{rough}}} \frac{1}{n^{1+iw}}$. But, since $h$ is multiplicative and $|h(p)| = 1$, we have (quite crudely)
$$ \Biggl|\sum_{\substack{n \leq x^{2\beta}, \\ n \; \text{is} \; P \; \text{rough}}} \frac{|h(n)|^2}{n^{1+iw}} - \sum_{\substack{n \leq x^{2\beta}, \\ n \; \text{is} \; P \; \text{rough}}} \frac{1}{n^{1+iw}} \Biggr| \ll \sum_{\substack{n \leq x^{2\beta}, \\ n \; \text{is} \; P \; \text{rough}, \\ n \; \text{is not squarefree}}} \frac{1}{n} \leq \sum_{p > P} \frac{1}{p^2} \sum_{\substack{n \leq x^{2\beta}/p^2 , \\ n \; \text{is} \; P \; \text{rough}}} \frac{1}{n} \ll \frac{\log x}{P\log P} , $$
and this tiny difference is negligible in all of the calculations. In the analogue of Key Propositions \ref{keyprop2} and \ref{keyprop4}, one needs an upper bound for $\E|\sum_{\substack{n \leq x^{2\beta}, \\ P \; \text{rough}}} \frac{h(n) f(n)}{n^{1/2}} |^2 |\sum_{\substack{n \leq x^{2\beta}, \\ P \; \text{rough}}} \frac{h(n) f(n)}{n^{1/2+iw}} |^2$ rather than $\E|\sum_{\substack{n \leq x^{2\beta}, \\ P \; \text{rough}}} \frac{f(n)}{n^{1/2}} |^2 |\sum_{\substack{n \leq x^{2\beta}, \\ P \; \text{rough}}} \frac{f(n)}{n^{1/2+iw}} |^2$. But again, examining the start of the proof of Probability Result \ref{probres3} we see the difference between these may be bounded by $\sum_{\substack{N \leq x^{4\beta}, \\ N \; \text{is} \; P \; \text{rough}, \\ N \; \text{is not squarefree}}} \frac{d(N)^2}{N} \ll \sum_{p > P} \frac{1}{p^2} \sum_{\substack{N \leq x^{4\beta}/p^2 , \\ N \; \text{is} \; P \; \text{rough}}} \frac{d(N)^2}{N} \ll \frac{1}{P} (\frac{\log x}{\log P})^4$, say, which is insignificant.

The only remaining concern is the handling of the $P$-smooth sum inside $F_{\chi, h}^{*}(s)$, and the accompanying barrier $G_{\chi, h}(s)$. The calculations from sections \ref{secpasstormf} and \ref{subsecgtow}, allowing one to replace $W_{f, h}(s)$ by $G_{f, h}(s)$ (now on the random multiplicative function side) where necessary, go through identically. Finally, since Lemmas 1 and 6 of Harper~\cite{harperrmflow} go through without change for the $h(n)$-weighted random Euler products (because only the first order terms $\frac{h(p) f(p)}{p^{s}}$ in the Taylor expansions of the logarithms contribute there, and our assumption that $|h(p)| = 1$ implies they contribute the same as $\frac{f(p)}{p^{s}}$), all of the arguments in sections 3.2 and 5.1 of Harper~\cite{harperrmflow} go through as well, and we can deduce the necessary substitutes for our Probability Results \ref{probres4} and \ref{probres5}.

Given the same bounds as in Key Propositions \ref{keyprop1} and \ref{keyprop2}, the argument with H\"older's inequality from the Introduction goes through identically to the unweighted case, completing the proof.
\qed

\subsection{The zeta sum case: Theorem \ref{thmlb3}}
To handle the zeta sum case, one can follow the proof of Theorem \ref{thmlb1} with $I(\chi)$ replaced by
$$ I(t) = I_{x,q}(t) := \frac{1}{2\pi i} \int_{1/2-i\mathcal{T}}^{1/2+i\mathcal{T}} \Biggl(\sum_{x^{1-3\beta} < p \leq x^{1-\beta}} \frac{p^{it}}{p^s} \Biggr) F_{t}^{*}(s) G_{t}(s) \frac{x^s}{s} ds , $$
where
$$ F_{t}^{*}(s) := \Biggl(\sum_{\substack{n \leq x^{\beta}, \\ n \; \text{is} \; P \; \text{smooth}}} \frac{n^{it}}{n^s} \Biggr) \cdot \Biggl(\sum_{\substack{n \leq x^{2\beta}, \\ n \; \text{is} \; P \; \text{rough}}} \frac{n^{it}}{n^s} \Biggr) , $$
and $G_{t}(s) = G_{t,q,P}(s) := \left( \prod_{0 \leq j \leq \log\log P - 1} \tilde{\gamma_{j}}\left( \Re \sum_{p \leq P^{e^{-j}}} \sum_{k=1,2} \frac{p^{itk}}{kp^{ks}} \right) \right)^{\lfloor \log\log P \rfloor}$ (as constructed in section \ref{secpasstormf}).

Unlike in the character sum setting, the ``continuous characters'' $n^{it}$ do not enjoy perfect orthogonality when integrated over $[0,T]$, instead we have the approximate orthogonality relation $\frac{1}{T} \int_{0}^{T} n^{it} m^{-it} dt = \frac{1}{T} \int_{0}^{T} e^{it\log(n/m)} dt = \textbf{1}_{n=m} + O(\frac{\textbf{1}_{n \neq m}}{T|\log(n/m)|})$. One could rework the character sum proofs with some additional (ultimately acceptable) error terms arising, but a cleaner approach is perhaps the following.

Let $\Phi : \R \rightarrow \R$ be the Beurling--Selberg function (see e.g. Vaaler’s paper~\cite{vaaler}) majorising the indicator function $\textbf{1}_{[0,1]}$, and with Fourier transform supported on $[-1/\delta, 1/\delta]$. We take $\delta = 1/\log T$, say. In particular, we then have $\frac{1}{T} \int_{-\infty}^{\infty} \Phi(t/T) dt = \int_{-\infty}^{\infty} \Phi(x) dx = 1 + \delta$, and for any $n,m \leq \delta T/2\pi$ with $n \neq m$ we get $\frac{1}{T} \int_{-\infty}^{\infty} n^{it} m^{-it} \Phi(t/T) dt = \frac{1}{T} \int_{-\infty}^{\infty} e^{-it\log(m/n)} \Phi(t/T) dt = \hat{\Phi}((T/2\pi)\log(m/n)) = 0$.  Then
$$ \frac{1}{T} \int_{0}^{T} |\sum_{n \leq x} n^{it}|^{2q} dt = \frac{1}{T} \int_{-\infty}^{\infty} \Phi(t/T) |\sum_{n \leq x} n^{it}|^{2q} dt - \frac{1}{T} \int_{-\infty}^{\infty} (\Phi(t/T) - \textbf{1}_{[0,T]}(t)) |\sum_{n \leq x} n^{it}|^{2q} dt , $$
where $(\Phi(t/T) - \textbf{1}_{[0,T]}(t))$ is always non-negative.

By H\"older's inequality and our observations above (and the fact that $x \leq T^{0.499}$ is certainly smaller than $\delta T/2\pi$), the subtracted term is
\begin{eqnarray}
& \leq & \Bigl(\frac{1}{T} \int_{-\infty}^{\infty} (\Phi(t/T) - \textbf{1}_{[0,T]}(t)) dt \Bigr)^{1-q} \Bigl(\frac{1}{T} \int_{-\infty}^{\infty} (\Phi(t/T) - \textbf{1}_{[0,T]}(t)) |\sum_{n \leq x} n^{it}|^{2} dt \Bigr)^{q} \nonumber \\
& = & \delta^{1-q} \cdot \Bigl(\sum_{n,m \leq x} \frac{1}{T} \int_{-\infty}^{\infty} n^{it} m^{-it} (\Phi(t/T) - \textbf{1}_{[0,T]}(t)) dt \Bigr)^{q} = \delta^{1-q} \Bigl(\sum_{n \leq x} \delta + O(\sum_{\substack{n,m \leq x, \\ n \neq m}} \frac{1}{T|\log(n/m)|}) \Bigr)^{q} . \nonumber
\end{eqnarray}
The second bracket is $\ll \delta x + \sum_{n < m \leq x} \frac{1}{T \log(m/n)} \ll \delta x + \sum_{n < m \leq x} \frac{m}{T (m-n)} \ll \delta x + \frac{x^2 \log(2x)}{T}$, which again is certainly $\ll \delta x$ when $\delta = 1/\log T$ and $x \leq T^{0.499}$. Thus the whole subtracted term is $\ll \delta x^q = \frac{x^q}{\log T}$, which is negligible, and we see it will suffice to prove Theorem \ref{thmlb3} for the weighted average $\frac{1}{T} \int_{-\infty}^{\infty} \Phi(t/T) |\sum_{n \leq x} n^{it}|^{2q} dt$.

But with the weight $\Phi(t/T)$ present, the functions $n^{it}$ {\em are} perfectly orthogonal for numbers $n \leq \delta T/2\pi = T/(2\pi \log T)$, as computed above. We can then pass to random multiplicative functions exactly as in section \ref{secpasstormf}, and deduce precise analogues of Key Propositions \ref{keyprop1} and \ref{keyprop2} for the $\Phi$-weighted averages of $\sum_{n \leq x} n^{it}$ and $I(t)$ (with the conditions $x^{1+2\beta+ \beta/(\log\log x)^{10}} < r$ and $x^{2+4\beta+ \beta/(\log\log x)^{10}} < r$ replaced by $x^{1+2\beta+ \beta/(\log\log x)^{10}} \leq T/(2\pi \log T)$ and $x^{2+4\beta+ \beta/(\log\log x)^{10}} \leq T/(2\pi \log T)$). Theorem \ref{thmlb3} follows from this, by exactly the same H\"older's inequality computations as in the Introduction.
\qed



\begin{thebibliography}{99}

\bibitem{bretechemunschten} R. de la Bret\`eche, M. Munsch, G. Tenenbaum. Small G\'al sums and applications. {\em J. Lond. Math. Soc. (2)}, \textbf{103}, no. 1, pp 336-352. 2021

\bibitem{gaowulow} P. Gao, X. Wu. Upper Bounds for low moments of twisted Fourier coefficients of modular forms. Preprint available online at \url{https://arxiv.org/abs/2512.05378}.

\bibitem{gorodetsky} O. Gorodetsky. Magic squares, the symmetric group and M\"{o}bius randomness. {\em Monatsh. Math.}, \textbf{204}, no. 1, pp 27-46. 2024

\bibitem{gransoundlcs} A. Granville, K. Soundararajan. Large character sums. {\em J. Amer. Math. Soc.}, \textbf{14}, no. 2, pp 365-397. 2001

\bibitem{gut} A. Gut. {\em Probability: A Graduate Course.} Second edition, published by Springer Texts in Statistics. 2013

\bibitem{shardy} S. Hardy. The distribution of partial sums of random multiplicative functions with a large prime factor. Preprint available online at \url{https://arxiv.org/abs/2503.06256}.

\bibitem{harperpartition} A. J. Harper. On the partition function of the Riemann zeta function, and the Fyodorov--Hiary--Keating conjecture. Preprint available online at \url{https://arxiv.org/abs/1906.05783}.

\bibitem{harperrmfhigh} A. J. Harper. Moments of random multiplicative functions, II: High moments. {\em Algebra Number Theory}, \textbf{13}, no. 10, pp 2277-2321. 2019

\bibitem{harperrmflow} A. J. Harper. Moments of random multiplicative functions, I: Low moments, better than squareroot cancellation, and critical multiplicative chaos. {\em Forum of Mathematics, Pi}, \textbf{8}, e1, 95pp. 2020

\bibitem{harpertypicalchar} A. J. Harper. The typical size of character and zeta sums is $o(\sqrt{x})$. Preprint available online at \url{https://arxiv.org/abs/2301.04390}.

\bibitem{harperas} A. J. Harper. Almost sure large fluctuations of random multiplicative functions. {\em Int. Math. Res. Not.}, \textbf{2023}, no. 3, pp 2095-2138.

\bibitem{harperrmf3} A. J. Harper. Moments of random multiplicative functions, III: A short review. In {\em Proceedings of the International Congress of Basic Science 2024}, International Press of Boston, Inc., also available online at \url{https://arxiv.org/abs/2410.11523}.

\bibitem{harperbtsqc} A. J. Harper. Better than squareroot cancellation in number theory. To appear in the Proceedings of the ICM 2026, available online at \url{https://arxiv.org/abs/2512.23681}

\bibitem{harperlowchar2} A. J. Harper. Lower bounds for low moments of character sums, II: Long unweighted sums, and non-vanishing. {\em In preparation}.

\bibitem{korner} T. W. K\"orner. {\em Fourier analysis.} Reprint, published by Cambridge University Press. 2022

\bibitem{mvchar} H. L. Montgomery, R. C. Vaughan. Mean values of character sums. {\em Canad. J. Math.}, \textbf{31}, no. 3, pp 476-487. 1979

\bibitem{mv} H. L. Montgomery, R. C. Vaughan. {\em Multiplicative Number Theory I: Classical Theory.} First edition, published by Cambridge University Press. 2007

\bibitem{szaboupper} B. Szab\'o. High moments of theta functions and character sums. {\em Mathematika}, \textbf{70}, no. 2, Paper No. e12242, 37 pp. 2024

\bibitem{szabolower} B. Szab\'o. A lower bound on high moments of character sums. Preprint available online at \url{http://www.arxiv.org/abs/2409.13436}.

\bibitem{tenenbaumintro} G. Tenenbaum. {\em Introduction to analytic and probabilistic number theory.} Third edition, translated from the 2008 French edition by Patrick D. F. Ion, published by the American Mathematical Society, Providence, RI. 2015

\bibitem{vaaler} J. Vaaler. Some extremal functions in Fourier analysis. {\em Bull. Amer. Math. Soc. (N.S.)}, \textbf{12}, no. 2, pp 183-216. 1985

\bibitem{wangxu} V. Y. Wang, M. W. Xu. Harper's beyond square-root conjecture. {\em Int. Math. Res. Not.}, \textbf{2025}, no. 18, 22 pp.



\end{thebibliography}
\end{document}